%% file: rand-grp-elts.tex
\documentclass[11pt]{article}
\usepackage{latexsym}  % To use \Box, compatibility with latex 2.09
\usepackage{epsfig}

\title{Towards a Practical, Theoretically Sound Algorithm for
                 Random Generation in Finite Groups}
\author{Gene Cooperman \\
         College of Computer Science\\
         Northeastern University\\
         Boston, MA 02115\\
         gene@ccs.neu.edu
}
\newcommand{\qed}{$\Box$}

\newtheorem{theorem}{Theorem}[section]
\newtheorem{corollary}[theorem]{Corollary}
\newtheorem{lemma}[theorem]{Lemma}
\newenvironment{proof}{\par\noindent{\it Proof:}\/\ \ }{ \qed\medskip}
\newtheorem{definition}{Definition}
\newtheorem{remark}{Remark}

\def\PrOP{\mathop{\rm Pr}}
\def\Pr(#1){\PrOP({#1})}
\def\E{\mathop{\rm E}\nolimits}
\def\supp(#1){\mathop{\rm supp}({#1})}
\def\GL{{\rm GL}}
\def\GF{{\rm GF}}
\def\Syl{{\rm Syl}}
\def\calP{{\cal P}}
\def\calR{{\cal R}}
\def\calT{{\cal T}}
\def\calX{{\cal X}}
\def\a{{a}} % variable a
\def\cbar{{\overline c}}
\def\mathbreak{$\linebreak[0]$}
\def\inv{^{-1}}
\def\gen#1{{\langle #1\rangle}}
\def\defeq{\mathrel{\mathop=\limits^{\rm def}}}
\def\probeq{\mathrel{\mathop=\limits^{\scriptscriptstyle\rm prob}}}
\def\norm#1{||#1||}

\long\def\IGNORE#1{}

\begin{document}

\maketitle

\begin{abstract}
\begin{sloppypar}
  This work presents a new, simple $O(\log^2|G|)$ algorithm, the
  Fibonacci cube algorithm, for producing random group elements in
  black box groups.  After the initial $O(\log^2|G|)$ group
  operations, $\varepsilon$-uniform random elements are produced using
  $O((\log 1/\varepsilon)\log|G|)$ operations each.  This is the first
  major advance over the ten year old result of Babai~\cite{Babai91},
  which had required $O(\log^5|G|)$ group operations.  Preliminary
  experimental results show the Fibonacci cube algorithm to be
  competitive with the product replacement algorithm.
\end{sloppypar}

  The new result leads to an amusing reversal of the state of affairs
  for permutation group algorithms.  In the past, the fastest random
  generation for permutation groups was achieved as an application of
  permutation group membership algorithms and used deep knowledge
  about permutation representations.  The new black box random
  generation algorithm is also valid for permutation groups, while
  using no knowledge that is specific to permutation
  representations.  As an application, we demonstrate a new algorithm
  for permutation group membership that is asymptotically faster than
  all previously known algorithms.
\end{abstract}

\section{Introduction}

Quickly finding an element of a black box group is a problem of
critical importance for many randomized algorithms for mathematical
groups.  (Black box groups are defined later.)  Random group elements
are especially important for computations with finite matrix groups,
where few efficient deterministic algorithms are known.

Researchers requiring generation of such random elements tended to
have a split personality.  On the one hand, one could chose a
theoretically sound algorithm with a complexity that was far too high
to be practical.  The best previous theoretical algorithm required
$O(\log^5|G|)$ group multiplications to produce a random
element~\cite{Babai91}.  On the other hand, one could choose a
heuristic for random elements such as the product replacement
algorithm~\cite{CellerLeedhamGreenEtAl95}, which could be demonstrated
to have a bias away from the uniform distribution~\cite{BabaiPak02},
but was ``good enough'' in practice.

This paper presents a simple $O(\log^2|G|)$ algorithm, the Fibonacci
cube algorithm, which is easy to program.  After the initial
$O(\log^2|G|)$ group operations, $\varepsilon$-uniform random elements
are produced using $O((\log 1/\varepsilon)\log|G|)$ operations each.
The algorithm is in Section~\ref{sec:pseudo-code}.
The main theoretical result is Theorem~\ref{thm:main}.  The
theoretical analysis of this paper currently has an unacceptably high
coefficient of complexity, although experimental results show it to be
competitive with the product replacement algorithm.  The conclusion
points out opportunities to lower the theoretical coefficient by
refining the complexity analysis.

A {\em black box group} is a group with an associated oracle, in which
group elements are encoded as binary strings of some uniform
length~$L$.  The oracle can multiply, find inverses, and compare an
element with the identity.  Note that this implies an upper bound
of~$2^L$ on the group order.

A common use of black box groups is to model finite matrix groups over
finite fields.  A matrix group, $\GL(d,q)$ (dimension~$d$ over
$\GF(q)$), is a black box group with an encoding of length
$L=d^2\log_2q$, and its order is a priori bounded by $2^L$.  Almost
every paper in the recent development of matrix group algorithms
assumes the availability of a random generation algorithm.  In
particular, the matrix recognition project~\cite{Leedham-Green01} is a
project to recognize matrix groups in $\GL(d,q)$ for values of $d$ up
to approximately~100, and for moderate size values of~$q$.  That
project relies heavily on the ability to compute nearly random group
elements.

Surprisingly, even in the regime of permutation groups, the new black
box algorithm for random generation is faster than the best know
permutation algorithm both for the case of large and small base.  Let
$n$ be the permutation degree.  For large base, $\log|G|\le n\log n$,
and so we have random generation in $O(n^3\log^2n)$.  For small base,
if we assume a base size of $O(\log n)$, then $\log|G|\le\log^2n$.

Let $G=\gen{S}$ be a finite black box group.  We use $\Pr(\cdot)$ to
notate probability and $\E(\cdot)$ to notate expectation.  Random
variables are denoted by capital letters, while group elements are
denoted by lower case letters.  Let $U$ be a random variable on~$G$
with uniform distribution.  We use the notation $A\subset G$ for a
{\em proper subset} of~$G$ and $A\subseteq G$ for a {\em subset}
of~$G$.  Similarly, $H<G$ denotes a {\em proper subgroup} of~$G$ and
$H\le G$ denotes a {\em subgroup} of~$G$.

% end IGNORE

\subsection{Previous work}

\IGNORE{

{\bf
The new algorithm reduces the complexity of random generation from
$O(\log^5|G|)$ to $O(\log^2|G|)$.  Specifically, it produces
$O(\log|G|)$ random group elements in $O(\log^2|G|)$.  Thereafter,
random elements are produced at a cost of $O(\log|G|)$ group operations.
The new algorithm has a simpler proof, not requiring the flows, eigenvalue
estimates, etc. from the Babai approach.

It should be emphasized that we produce truly uniform random group elements
(rather than the $\varepsilon$-uniform elements of the Babai algorithm),
and we do it with exponential likelihood.  Hence, we upgrade the algorithm
from Monte Carlo to exponentially likely.  The probability that
we make an error in the preprocessing is $O(\log^{-2}|G|)$ since we just
need to force the maximum probability density, $M$, down to {\bf WHAT?}.
MUST ADD STOPPING CONDITION TO THEORY.  After that an additional
$\log^2|G|$ goes toward extending Chernoff bound.
For Babai, probability of error is $O(1/\log|G|)$.
Our probability of error should be $O(\exp(-\log^2|G|))$ or something
(exponentially small in time spent),
by Chernoff's bound.
}

} % end IGNORE

\begin{sloppypar}
The first polynomial time algorithm for random group
elements was demonstrated by Babai~\cite{Babai91}.  It runs in
time~$O(\log^5|G|)$.
Unfortunately, the high complexity means that this algorithm is not
used in computations.  As Babai wrote
in the {\sl Handbook of Combinatorics}~\cite{Babai95}:
\begin{quote}
Reducing the exponent~$5$ would be of great significance since
many algorithms in computational group theory rely on ``randomly
chosen'' elements from the group.~\cite{Babai95}
\end{quote}
\end{sloppypar}

A second heuristic, {\em product replacement}, was then proposed by
Celler et al.~\cite{CellerLeedhamGreenEtAl95} as a practical way to
find random elements of~$G$.

\IGNORE{

{\bf The original authors noted a slight bias, but
the complexity to approach the limiting distribution is not known?}

} % end IGNORE

Other researchers asked how fast a product replacement algorithm would
approach a uniform distribution in the class of generating $k$-sets
for~$G$.  Note that such a random $k$-set is distinct from a random
group element.  Diaconis and Saloff-Coste showed the algorithm to
produce nearly random generating $k$-sets in sub-exponential
time~\cite{DiaconisSaloffCoste98}, and Pak then showed it to operate
in polynomial time~\cite{Pak00}.  Pak requires the use of a $k$-tuple
in which $k=\Omega(\log|G|\log\log|G|)$.  When
$k=\theta(\log|G|\log\log|G|)$, he achieves his best time of
$O(\log^9|G|(\log\log|G|)^5)$.  Babai and Pak~\cite{BabaiPak02}
presented an important obstacle, whereby for the limiting
distribution of $k$-sets, individual group elements are shown to be
biased away from the identity.

\IGNORE{

  (Note that this is different from the problem of finding
random generating $k$-sets.)  In this case, they produce an infinite
family of groups that exhibits large bias.  As an example, they
propose the direct product group $G=(A_5)^{19}$.  $G$ can be generated
by two elements, and product replacement with 2-sets should exhibit
strong biases.  Intuitively, the product replacement algorithm
exhibits a small bias on~$A_5$, and the 19~components of the direct
product magnify the bias up to 19~times.  Nevertheless, in defense of
current heuristic practice, applications of the product replacement
algorithm typically use $k$-sets with $k\gg 2$.

} % end IGNORE

\subsection{Outline of contents}
\label{sec:outline}

The primary result is Theorem~\ref{thm:main}.
Informally, it shows that one can construct an initial, nearly random element
using $O(\log^2|G|)$ group operations with further elements produced
in $O(\log|G|)$ time.  The algorithm is given in
Section~\ref{sec:pseudo-code}.  For a high level overview of the
approach to the proof, see Section~\ref{sec:overview-proof} after
reading this section.

We use the notation $XY$ for the product of two $G$-valued random
variables in analogy with $gh$ for $g,h\in G$.  To illustrate the
notation, if $E$ is a $\{0,1\}$-random variable, then $XY^E=h$ if
$X=h$ and $E=0$, while $XY^E=hg$ if $X=h$, $Y=g$ and $E=1$.

Although the algorithm for generating random group elements is a
simple one, its justification is not simple.  We will develop a
sequence of random variables $\calR_0,\calR_1,\ldots$ such that
$\calR_0$ is fixed at the identity and
$\calR_{i+1}=\calR_i g_i^{E_i}$ where $g_i\in G$ is chosen from a
random distribution based on $\calR_i$ and $E_i$ is a uniform random
variable on $\{0,1\}$.  The $E_1,E_2,\ldots$ are pairwise independent
and independent of the other random variables.
For some fixed~$t=\Omega(\log|G|)$, $g_1^{E_1}\cdots g_t^{E_t}$ is a nearly
uniform random variable on~$G$, and computing an element from its
distribution requires at least~$t$ group operations.

Section~\ref{sec:preliminaries} provides some easy, well-known lemmas
which form the foundation for the rest of the paper.  The $\ell^2$
norm, $\norm{X}$, of a random variable~$X$ is defined in
Section~\ref{sec:norm-ell2}, along with some easy lemmas about it.
One easily shows that $\norm{\calR_i}$ is monotonically non-increasing
as a function of~$i$.  The $\ell^2$ norm had previously been used by Diaconis
and Saloff-Coste to analyze random walks on
groups~\cite{DiaconisSaloffCoste93}.

The primary goal of the proof is to show that $\norm{\calR_{i+1}}\le
c\norm{\calR_i}$ with probability at least $\rho>0$ for some positive
$c<1$.  Section~\ref{sec:overview-proof} outlines the ideas of that
proof.  Section~\ref{sec:proof} provides that proof
and concludes with the formal statement in
Lemma~\ref{lem:main} showing that for $t=\Omega((1/c)\log|G|)$,
$\calR_t$ is semi-uniform.  ($\Pr(\calR_t=g)$ is bounded away from~0.)  That
lemma then yields the main theorem, Theorem~\ref{thm:main}.

The main result relies on some
technical results from Sections~\ref{sec:reduction}
to~\ref{sec:escape}.  As a matter of notation, we reserve upper case
letters $E$, $I$, $J$, $K$, $T$ through~$Z$ and $\calR_i$ for random variables.

Section~\ref{sec:reduction} makes various assumptions on $X$ and $W$.
It then asks for what positive $\alpha<1$ and $\rho>0$ can one conclude
that for $g$ drawn from the distribution of~$W$,
$\norm{Xg^E}<\alpha\norm{X}$ with probability at least~$\rho$.
Here, $E$ is an independent random variable on~$\{0,1\}$.

Section~\ref{sec:decomposition} asks the following question.  Let $X$
and $Y$ be random variables on~$G$ and let $J$ be a random variable on
$\{0,1\}$.  (Often, we will take $\Pr(J=0)=1/2$ and $\Pr(J=1)=1/2$.)
Assume that $\calR_i$ has the same probability density function as
$X^JY^{1-J}$.  (Note that the notation $X^J$ means that
$X^J=X$ when $J=1$, and $X^J$ is the group
identity element when $J=0$.)
If $\norm{X g_i^{E_i}}\le\alpha\norm{X}$ for some
$0<\alpha\le 1$, then for what $\beta$ is it true that $\norm{\calR_i
  g_i^{E_i}}\le\beta\norm{\calR_i}$?  In order to state the results
more generally, that section writes $Z$ for $\calR_i$ and $W$ for
$g_i^{E_i}$.  As will be seen, Section~\ref{sec:algorithm}.

% {\bf
% The first idea, described in Section~\ref{sec:decomposition}, is that
% given a random variable $\calR_i$ on~$G$, we can usually choose a new
% random variable~$W$ on~$G$ such that $\calR_i$ decomposes into an
% identically distributed random variable, $W^KT^{1-K}$, for $T$ another
% random variable on~$G$ and for $K$ some
% $\{0,1\}$-random variable.  Hence, given $g\in G$ drawn from the
% distribution of~$\calR_i$, with probability $\Pr(K=1)$ we can treat it
% as coming from the distribution of~$W$.
% }

Section~\ref{sec:escape} worries about the unusual case of being
``stuck'' in a proper subgroup.  The
fixed $\{g_1,\ldots,g_i\}$ constructed to define the series
$\calR_1,\ldots,\calR_{i+1}$ can all be contained in a proper
subgroup.  In such a case, a random $g_{i+1}$ drawn from $\calR_{i+1}$
will also be in the proper subgroup.  Further, if $\calR_{i+1}$ has
only a small probability of lying outside a proper subgroup, then the
same problem arises.  The solution is to use the generators
of~$G$ to construct a group element~$g_i\notin A$, whereupon $X
g_i^{E_i}$ is smaller than~$X$ in the $\ell^2$ norm.  We use random
subproducts (Definition~\ref{def:randomSubproduct}) as an efficient
way to  construct a $g_i\notin A$.

Section~\ref{sec:escape} contains Theorem~\ref{thm:escape}, which may
have independent interest.  Informally, it states that for a set
$A\subseteq G$ with $A=A\inv$, either a random $(u,v)\in A\times A$
satisfies $uv\notin A$ with at least some positive probability, or
else $A$ is close to a subgroup~$A'$ with $A'A'$ a subgroup of~$G$.

Section~\ref{sec:algorithm} demonstrates the Fibonacci Cube algorithm,
which constructs the $g_i$ in the definition of $\calR_i$.  This is
the main algorithm.  This is enough to show that
$\calR_t\inv\overline\calR_t$ is semi-uniform for sufficiently
large~$t$.

Section~\ref{sec:semi-to-uniform} shows how to construct
$\varepsilon$-uniform random elements from $\varepsilon$-semi-uniform
random elements.  It then summarizes the previous results in the main
theorem, Theorem~\ref{thm:main}.  Theorem~\ref{thm:semiuniform} of
that section is of
independent interest, since it shows how to efficiently construct a
uniform random variable from a semi-uniform random variable.

Section~\ref{sec:experiments} presents some initial experimental
results applying the Fibonacci cube algorithm to conjugacy classes.
After a precomputation of about 100 group operations, one
produces independent pseudo-random elements costing 20 group
operations per random element.  Those elements satisfy the $\chi^2$
goodness of fit test as having a distribution over the conjugacy
classes that is close to uniformly random.

Section~\ref{sec:productReplacement} produces a $O(\log^2|G|)$ random
generation algorithm for a variation of the product replacement
algorithm, and Section~\ref{sec:permGroupMembership} describes how to
use the new Fibonacci cube algorithm to produce what is currently
asymptotically fastest group membership algorithm --- both for the
general (large base) case and the special case of small base groups.

\section{Preliminaries}
\label{sec:preliminaries}

The following easy lemmas and theorems are included for completeness.
Note that throughout this paper, random variables are always denoted
by upper case letters $E$, $I$, $J$, $K$, $T$ through~$Z$ and by~$\calR_i$.

\subsection{Probability and $\varepsilon$-uniform random variables}

The following lemma is well-known and has an easy proof.
\begin{lemma}[Markov's inequality]
\label{lem:markov}
Let $\xi$ be a nonnegative random variable and $\lambda>1$ a real number.
Then
\[
\Pr(\xi\ge\lambda\E(\xi))\le{1\over\lambda}.
\]
\end{lemma}

\begin{corollary}
\label{cor:markovAlt}
Let $\xi$ be a random variable on the interval~$[0,1]$ and $\lambda>1$
a real number.  Then
\[
\Pr(\xi>1-\lambda\E(1-\xi))\ge1-{1\over\lambda}.
\]
\end{corollary}
\begin{proof}
Let $\zeta=1-\xi$ and note that $\zeta$ is nonnegative.  Then
$
 \Pr(\xi>1-\lambda\E(1-\xi))\ge1-{1\over\lambda}\Leftrightarrow
 \Pr(1-\zeta>1-\lambda\E(\zeta))\ge1-{1\over\lambda}\Leftrightarrow
 \Pr(\zeta<\lambda\E(\zeta))\ge1-{1\over\lambda}\Leftrightarrow
 \Pr(\zeta\ge\lambda\E(\zeta))\le{1\over\lambda}
$,
and the last inequality follows from Markov's inequality.
\end{proof}

\begin{theorem}[Chernoff's Bound~\cite{Chernoff52}]
\label{thm:chernoff}
Let $S_t$ be a random variable
equal to the number of successes in $t$ independent Bernoulli trials
in which the probability of success is $p$ $(0 < p < 1)$. Let $0 <
\epsilon < 1$. Then
\[
\Pr(S_t \le \lfloor (1-\epsilon)pt \rfloor) \le e^{-{\epsilon}^2 pt/2}.
\]
\end{theorem}

\begin{definition}
\label{def:randomSubproduct}
A {\em random subproduct} on an ordered
set~$S=\{g_1,\ldots,g_k\}\subseteq G$ is given by
$g_1^{\epsilon_1}\cdots g_k^{\epsilon_k}$ for $\epsilon_i$
independent, uniform random variables on~$\{0,1\}$.
$(\Pr(\epsilon_i=0)=1/2$ and $\Pr(\epsilon_i=1)=1/2$.$)$
\end{definition}

\IGNORE{

\begin{quotation}
\noindent
========================================
\par\noindent
{\bf The issue is that maybe $Au\cap A=\emptyset$, $Ag_j\cap A=\emptyset$,
and $Aug_j\approx Au\approx Av\inv$.
If we can do this, then take a random subproduct of $1/\varepsilon$ 
elements, $uv$, and of each generator.
}
\begin{lemma}[random subproduct]
Let $A\subseteq G$ and let $S\subseteq G$ be a set having at least one
element~$s$ such that $As\cap A=\emptyset$.  Let $r$ be a random
subproduct on~$S$.  Then with probability at least 1/2, $|Ar\setminus
A|\ge|A|/2$.
\end{lemma}
\begin{proof}
Let $S=\{g_1,\ldots,g_k\}$ and let $j\le k$ be the largest integer
such that $Ag_j\cap A=\emptyset$.  Decompose the random subproduct
$r=g_1^{\epsilon_1}\cdots g_k^{\epsilon_k}$ as $r=ug_j^{\epsilon_j}v$.
If $|Aug_j\setminus Av\inv|\ge|A|/2$, then with probability 1/2,
$\epsilon_j=1$,
which implies $|Ar\setminus A|=|Aug_j\setminus Av\inv|\ge|A|/2$.

If $|Aug_j\setminus Av\inv|<|A|/2$, then with probability 1/2,
$\epsilon_j=0$.  Either $|Au\setminus Av\inv|<|A|/2$ (contradiction??)
or $|Au\setminus Av\inv|\ge|A|/2$ (then easy proof follows)
{\bf How does $Au$ compare with $A$?  Did it move much?
If it did, can we stop here, and take $r=u$?}
\end{proof}

{\it
If $|Aug_j\setminus Av\inv|<|A|/2$, then with probability 1/2,
$\epsilon_j=0$, which implies
$|Ar\setminus A|=|Auv\setminus A|=|Au\setminus Av\inv|
\ge |Au|-|Au\cap Av\inv|>|A|/2$.
}
{\bf Then next lemma is a corollary of this.}
\par\noindent
========================================
\end{quotation}

} % end IGNORE

The following is a generalization of Proposition~2.1 of Cooperman and
Finkelstein~\cite{CoopermanFinkelstein93}.
\IGNORE{
{\bf IS THIS ALSO THE LEMMA USED BY BLS FOR THEIR FOCS GROUP MEMBERSHIP?}
} % end IGNORE

% CWI: \cite{CoopermanFinkelstein92}

\begin{lemma}[random subproduct]
\label{lem:randomSubproduct}
Let $H$ be a proper subgroup of $G=\gen{S}$ and let $r$ be a random subproduct
on~$S$.  Then with probability at least 1/2, $|Hr\setminus H|\ge|H|/2$.
\end{lemma}
\begin{proof}
Let $S=\{g_1,\ldots,g_k\}$ and let $j\le k$ be the largest integer
such that $g_j\notin H$.  Decompose the random subproduct
$r=g_1^{\epsilon_1}\cdots g_k^{\epsilon_k}$ as $r=ug_j^{\epsilon_j}v$.
If $|Hu\setminus H|\ge|H|/2$, then with probability 1/2, $\epsilon_j=0$,
which implies $|Hr\setminus H|=|Hu\setminus H|$.  If $|Hu\setminus
H|<|H|/2$, then with probability 1/2, $\epsilon_j=1$, which implies
$|Hr\setminus H|=|Hug_j\setminus H|\ge |Hu\cap H|=|H|-|Hu\setminus H|>|H|/2$.
\end{proof}

\IGNORE{

{\bf
Decomposition 1:  Choose $U$ to be the peak and $V$ to be
everything else.  Prove that the peak gets moved and so the maximum
probability density is reduced by some fraction.
}

{\bf
Decomposition 2:
Choose $U$ to be the uniform distribution (perhaps with probability
density $1/(2|G|)$) and $V$ to be the left over part.  With constant
probability, $X$ chooses a value from~$U$.  After that, the product of
a uniform random variable and any other distribution is still uniform.
So, we need only take $k$ elements from $X_1$, \ldots, $X_k$ and hope
that $I_j=1$ for some $1\le j\le k$.
}

} % end IGNORE

\begin{sloppypar}
\begin{lemma}
\label{lem:maxProb}
Let $X$ and $Y$ be independent random variables on~$G$.  Then $\min_{h\in
G}\Pr(X=h)\le\Pr(XY=g)\le\max_{h\in G}\Pr(X=h)$ for all $g\in G$.
Similarly, $\min_{h\in
G}\Pr(Y=h)\le\Pr(XY=g)\le\max_{h\in G}\Pr(Y=h)$.
\end{lemma}
\begin{proof}
Note $\Pr(XY=g)=\sum_{h\in G}\Pr(X=h)\Pr(XY=g)=\sum_{h\in G}\Pr(X=h\mbox{ and
}Y=h\inv g)=\sum_{h\in G}\Pr(X=h)\*\Pr(Y=h\inv g)$.  Further,
$\min_{h\in G}\Pr(Y=h)=\min_{h\in G}\Pr(X=h)\*\sum_{f\in G}\Pr(Y=f\inv
g)\le\sum_{f\in G}\Pr(X=f)\*\Pr(Y=f\inv g) \le\max_{h\in
G}\Pr(X=h)\*\sum_{f\in G}\Pr(Y=f\inv g)=\max_{h\in G}\Pr(X=h)$.
A similar argument holds for $\min_{h\in G}\Pr(Y=h)$ and $\max_{h\in
G}\Pr(Y=h)$.
\end{proof}
\end{sloppypar}

\begin{lemma}[Babai and Szemer\'edi~\cite{BabaiSzemeredi84}]
\label{lem:cubeDoubling}
The following holds:  $g\notin A\inv A\Leftrightarrow {Ag\cap A}=\emptyset
        \Leftrightarrow |Ag\setminus A|=2|A|$.
\end{lemma}
The proof is clear.

\begin{definition}
A random variable on a group~$G$ is an $\varepsilon$-uniform random
variable if $|\Pr(X=g)-1/|G|\,|\le\varepsilon/|G|$ for all $g\in G$.
Note that a $0$-uniform random variable is just a uniform random variable.
\end{definition}

\begin{lemma}[$\varepsilon$-uniform random variable]
\label{lem:uniform}
Let $U$ and $V$ be independent random variables on a group~$G$ and let
$\varepsilon\ge0$.  If
$U$ is an $\varepsilon$-uniform random variable, then
$UV$ and $VU$ are also $\varepsilon$-uniform.
\end{lemma}
\begin{proof}
$|\Pr(UV=g)-1/|G||=\sum_{h\in G}(|\Pr(U=h)-1/|G|)\,\*\Pr(V=h\inv g)
=(|\Pr(U=h)-1/|G|)(\sum_{h\in G}\Pr(V=h\inv g))\le\varepsilon 1/|G|$.
A similar argument follows for $VU$.
\end{proof}

The next lemma shows that once a random variable~$U$ is found to be
uniform on~$A$ for $|A|>|G|/2$, $U\inv VU$ is $\varepsilon$-uniform
for arbitrary random variable~$V$.

\begin{lemma}
\label{lem:cubeUniform}
Let $\alpha$ be a constant satisfying $1/2<\alpha\le 1$.
Let $A$ be a subset of a group~$G$ such that $|A|\ge\alpha|G|$.  Let
$U_1$, $U_2$ and~$V$ be independent random variables on~$G$.
Let $U_1$ and $U_2$ be uniform on~$A$
with $\Pr(U_1=g)=\Pr(U_2=g)=0$ for $g\notin A$.
Then
\[\forall g\in G,\quad
  \frac{1-\alpha}{\alpha} \frac{1}{|G|}
  \ge\Pr(U_1\inv VU_2=g)-\frac{1}{|G|}
  \ge - \left(\frac{1-\alpha}{\alpha}\right)^2 \frac{1}{|G|}.
\]
Hence, $U_1\inv VU_2$ is a $(1-\alpha)/\alpha$-uniform random variable on~$G$.
\end{lemma}
\begin{proof}
Note that $|A|\ge\alpha|G|$ implies $|A\cap Ag|\ge({2\alpha-1})|G|$.  So,
$\Pr(VU_2\in Ag)\ge({2\alpha-1})\*|G|/|A|\ge({2\alpha-1})/\alpha$.  Since $U_1$
and $U_2$ are independent, $\Pr(U_1\inv
VU_2=g)=\Pr(VU_2=U_1g)=\Pr(VU_2\in Ag)/|A|
\ge\left(({2\alpha-1})/\alpha^2\right)/|G|$.
Also $\Pr(U_1\inv VU_2=g)=\Pr(VU_2\in Ag)/|A|\le1/|A|=(1/\alpha)/|G|$.
Subtracting $1/|G|$ from the lower and upper bounds on $\Pr(U_1\inv
VU_2=g)$ completes the proof.
\end{proof}

In fact, Lemma~\ref{lem:cubeUniform} can easily be generalized to
$U_1$ uniform on~$A_1$ for $|A_1|\ge\alpha_1|G|$ and
$U_2$ uniform on~$A_2$ for $|A_2|\ge\alpha_2|G|$, but the existing
form suffices for our purposes.

\IGNORE{

The next two easy theorems are included for completeness.

\begin{theorem}
\label{thm:boundedMax}
Let $X$ and $Y$ be independent random variables on a finite group~$G$.
Let $m=\max_{g\in G}\Pr(X=g)$.  Then  $\Pr(XY=g)\le m$
for all $g\in G$.
\end{theorem}
\begin{proof}
For all $g\in G$, $\Pr(XY=g)=\sum_{h\in
G}\Pr(X=gh\inv)\Pr(Y=h)\le\sum_{h\in G}m\Pr(Y=h)=m$.
\end{proof}

} % end IGNORE

\subsection{The $\ell^2$ norm}
\label{sec:norm-ell2}

Let $R$ denote the real numbers.
Recall that the $\ell^2$ norm on $v=(v_1,\ldots,v_k)\in R^k$
is $\norm{v}_2=\sqrt{\sum_{i=1}^k (v_i)^2}$.
Let $\calX$ be the set of $G$-valued random variables for $G$~a group.
Define the function $\varphi$ as the natural function
from $\calX$ to~$R^{|G|}$,
the $|G|$-dimensional vector space over the reals.
Hence, if $X\in\calX$ and $G=\{g_1,g_2,\ldots,g_{|G|}\}$,
then define:
\begin{eqnarray*}
\varphi(X)&=&\left(\Pr(X=g_1),\Pr(X=g_2),\ldots,\Pr(X=g_{|G|})\right)\\
\norm{X}&=&\norm{\varphi(X)}_2=\sqrt{\sum_{g\in G}\left(\Pr(X=g)\right)^2}
\end{eqnarray*}
Note that $\norm{XY}$ is a norm under multiplication, since
$\norm{XY}=\norm{\varphi(XY)}_2
\le\norm{\varphi(X)}_2\,\norm{\varphi(Y)}_2=\norm{X}\,\norm{Y}$ by the
Cauchy-Schwartz inequality.

Observe that for two $G$-valued random variables $X$ and~$Y$,
\[
\norm{XY}=\norm{\sum_{g\in G}\varphi(Xg)\Pr(Y=g)}_2
        =\norm{\sum_{g\in G}\varphi(gY)\Pr(X=g)}_2.
\]

\begin{lemma}
\label{lem:equivalentNorms}
For $X$ a random variable on the group~$G$ and $g\in G$,
$\norm{X}=\norm{X\inv}=\norm{Xg}$.
\end{lemma}
\begin{sloppypar}
\begin{proof}
$\norm{X}^2=\sum_{h\in G}(\Pr(X=h))^2 =\sum_{h\in
G}(\Pr(X\inv=h\inv))^2=\norm{X\inv}^2$.  Similarly,
$\norm{X}^2=$\linebreak[0]$\sum_{h\in G}(\Pr(X=h))^2 =\sum_{h\in
G}(\Pr(Xg=hg))^2=\norm{Xg}$.
\end{proof}
\end{sloppypar}

\begin{lemma}
\label{lem:bounded-ell2}
If $X$ and $Y$ are independent $G$-valued random variables for
$G$ a group, then
$\norm{XY}\le\min(\norm{X},\norm{Y})$.
\end{lemma}
\begin{proof}
By the triangle inequality,
$\norm{XY}=\norm{\sum_{g\in G}\varphi(Xg)\Pr(Y=g)}_2
        \le\sum_{g\in G}\norm{\varphi(Xg)}_2\*\Pr(Y=g)
        =\norm{X}\*\sum_{g\in G}\Pr(Y=g)
        =\norm{X}$,
and similarly $\norm{XY}\le\norm{Y}$.
\end{proof}

\begin{lemma}
\label{lem:maxNorm-ell2}
Let $X$ be a random variable on~$G$.
If $\Pr(X=g)\le m$ for all $g\in G$, then $\norm{X}\le\sqrt{m}$.
\end{lemma}
\begin{proof}
$\norm{X}$ is maximized when $\Pr(X=g)=m$ or $\Pr(X=g)=0$ for all
$g\in G$ except at most one $g'\in G$ for which $0<\Pr(X=g')<m$.
To see this, let $Y$ be a random variable with $\Pr(Y=g)\le m$ such
that $\norm{Y}$ is maximal.
If $x_1=\Pr(Y=g_1)$, $x_2=\Pr(Y=g_2)$, $0<x_1<m$, $0<x_2<m$ and
$0<\delta\le x_2$, then $(x_1+\delta)^2 + (x_2-\delta)^2
=x_1^2+x_2^2+2(x_1-x_2)\delta+2\delta^2>x_1^2+x_2^2$ when $x_1>x_2$.
This violates maximality of $\norm{Y}$.  So there is at most one $g'\in G$ such
that $0<\Pr(Y=g')<m$.  Let $\Pr(Y=g')=m'<m$.  Then $\norm{Y}=\sqrt{{m'}^2
+ ((1-m')/m)m^2}<\sqrt{m}$.
\end{proof}

\begin{definition}
The {\em support} of a random variable~$X$ on a group~$G$ is the set
\[
\supp(X)=\{g\in G\colon~\Pr(X=g)>0\}.
\]
\end{definition}

\begin{lemma}
\label{lem:minNorm-ell2}
Let $X$ be a random variable on~$G$.
Then $\norm{X}\ge 1/\sqrt{|\supp(X)|}$.
\end{lemma}
\begin{proof}
Let $U$ be the uniform random variable on $\supp(X)$ and observe that
$\Pr(U=g)=1/|\supp(X)|$ for $g\in\supp(X)$.  Taking the inner product of
$\phi(U)$ and $\phi(X)$, the result follows from
$1/|\supp(X)|=\phi(U)\cdot\phi(X)\le\norm{U}\,\norm{X}
=\norm{X}/\sqrt{|\supp(X)|}$.
The inequality $\phi(U)\cdot\phi(X)\le\norm{U}\,\norm{X}$ is the
Cauchy-Schwartz inequality.
\end{proof}

\section{Reduction of probability in the $\ell^2$ norm}
\label{sec:reduction}

In this section, we derive estimates of the form
$\norm{Xg^{E}}\le\alpha\norm{X}$ for 
 ${E}$ a uniform $\{0,1\}$-random variable and for fixed $g$ drawn from the
distribution of~$W$, with probability at least $\rho>0$.
The positive parameters $\alpha<1$ and~$\rho$ depend on the choice
of~$X$ and~$W$.  In applications, we will find $X$, $Y$ and~$J$ such
that $\calR_i$ has the same distribution as $X^JY^{1-J}$.  Having
shown $\norm{Xg^{E}}\le\alpha\norm{X}$ in this section, 
Section~\ref{sec:decomposition} will allow us to conclude
$\norm{\calR_i g_i^{E_i}} < \alpha\norm{\calR_i}$ with probability at
least $\rho>0$.

\begin{lemma}
\label{lem:expectedNormReduction}
Let $X$ and $W$ be independent random variables on a group~$G$.  Let $E$ be a
$\{0,1\}$-random variable and let $X$ and $E$ be independent.  The
notation $\E_{g\in W}(f(g))$ denotes $\E(f(W))$ for the function
$f\colon~G\rightarrow R$ into the real numbers~$R$.  Hence,
$\E_{g\in W}(\norm{Xg^E}^2) \defeq \E(f(W))$ for
$f(g)=\norm{Xg^E}$ for $g\in G$.  Then
\begin{eqnarray*}
\E_{g\in W}(\norm{Xg^E}^2)
        &=& \left(\bigl(\Pr(E=0)\bigr)^2
                                + \bigl(\Pr(E=1)\bigr)^2\right) \norm{X}^2 + \\
        &&  \sum_{h\in G} 2 \Pr(E=0) \Pr(E=1) \Pr(X=h) \Pr(XW=h)
\end{eqnarray*}
\end{lemma}
\begin{proof}
Lemma~\ref{lem:equivalentNorms} tells us that
$\norm{X}=\norm{Xg\inv}$.
Without loss of generality, we can take $X$ and $W$ as independent.
If $X$ and $W$ were dependent, then
we would take $X'$ as an independent random variable with identical
distribution to~$X$, and note that
$\E_{g\in W}(\norm{X'g^E}^2)=\E_{g\in W}(\norm{Xg^E}^2)$.
For $X$ and $W$ independent,
$\sum_{g\in G}\Pr(X=hg\inv)\Pr(W=g)=\Pr(XW=h)$.
The following equality then holds.
\begin{eqnarray*}
\E{g\in W}(\norm{Xg^E}^2)
        &=& \sum_{g\in G} \left(\Pr(W=g)
          \sum_{h\in G} \left(\Pr(Xg^E=h)\right)^2\right)\\
%=& \sum_{g\in G} \left(\Pr(W=g)
%   \sum_{h\in G} \left(\Pr(X=h)\Pr(E=0)
%    + \Pr(X=hg\inv)\Pr(E=1)\right)^2\right)\\
        &=& \sum_{g\in G}\sum_{h\in G} \Pr(W=g) \biggl(
               \Bigl(\Pr(E=0)\Pr(X=h)\Bigr)^2
                 + \left(\Pr(E=1)\Pr(X=hg\inv)\right)^2 +\\
        &&  \qquad\qquad\qquad\qquad \quad
                 2 \Pr(E=0) \Pr(E=1) \Pr(X=h) \Pr(X=hg\inv) \biggr)\\
        &=& \sum_{g\in G} \Pr(W=g) \bigl(\Pr(E=0)\bigr)^2 \norm{X}^2
          + \sum_{g\in G} \Pr(W=g) \bigl(\Pr(E=1)\bigr)^2 \norm{Xg\inv}^2 +\\
        &&  \sum_{h\in G} 2 \Pr(E=0) \Pr(E=1) \Pr(X=h) \Pr(XW=h)\\
        &=& \left(\bigl(\Pr(E=0)\bigr)^2
                                + \bigl(\Pr(E=1)\bigr)^2\right) \norm{X}^2 +\\
        &&  \sum_{h\in G} 2 \Pr(E=0) \Pr(E=1) \Pr(X=h) \Pr(XW=h)
\end{eqnarray*}
\end{proof}

\begin{theorem}
\label{thm:specificReductionBySupport}
Let $X$, $W$ and~$Z$ be random variables on a group~$G$.  Let $E$ be a
uniform $\{0,1\}$-random variable and let $X$ and $E$ be independent.
Let $\lambda>1$ and let $g\in G$ be drawn from the
distribution of~$W$.  Let $\phi=\Pr(XW\in{\supp(X)})$.  Let $Z$ have
a density function such that $\Pr(Z=g)=\Pr(XW=g)/\phi$ for $g\in\supp(X)$
and $\Pr(Z=g)=0$ for $g\notin\supp(X)$.
(The random variable~$Z$ can be thought of as $XW$ conditioned on the
event $XW\in\supp(X)$.)
Let $\norm{Z}\le c\norm{X}$.
Then with probability at least
$1-1/\lambda$,
\[
        \norm{Xg^E}<\sqrt{\lambda\frac{1+c\phi}{2}}\,\norm{X}.
\]
\end{theorem}
\begin{proof}
Note $\sum_{h\in G}\Pr(X=h)\*\Pr(XW=h)
\le \sum_{h\in G}\Pr(X=h)\*\Pr(Z=h)\phi
\le \phi\norm{X}\,\norm{Z}\le c\phi\norm{X}^2$, where the first inequality
holds due to the Cauchy-Schwartz inequality.
From Lemma~\ref{lem:expectedNormReduction},
$\E_{g\in W}(\norm{Xg^E}^2)\le\norm{X}^2(1+c\phi)/2$.
Define the function $f(g)=\norm{Xg^E}^2$ from $G$ to the real numbers.
By Markov's inequality, $\PrOP(f(W)\ge
\lambda(\norm{X}^2(1+c\phi)/2))\le 1/\lambda$, from which the theorem
follows.
\end{proof}

The estimate of the next corollary is used for Case~2 in
Section~\ref{sec:algorithm}.

\begin{corollary}
\label{cor:altReductionBySupport}
\begin{sloppypar}
Assume the same hypotheses as Theorem~\ref{thm:specificReductionBySupport},
with the exception that $XW$ is replaced by $WX$ in the definition
of~$\phi$ and of~$Z$.
Then with
probability at least $1-1/\lambda$,
\[
        \norm{g^EX}<\sqrt{\lambda\frac{1+c\phi}{2}}\,\norm{X}.
\]
\end{sloppypar}
\end{corollary}
\begin{sloppypar}
\begin{proof}
  Replace $X$ by $X\inv$, $W$ by $W\inv$, and $g$ by $g\inv$ in
  Theorem~\ref{thm:specificReductionBySupport}.  Then
  $\phi=\Pr(WX\in{\supp(X)})$ and $\Pr(Z\inv=g)=\Pr(WX=g)/\phi$ for
  $g\in\supp(X)$ and $\Pr(Z\inv=g)=0$ otherwise.  So
  $\norm{Z\inv}=\norm{Z}\le c\norm{X}$.  Also
  $\norm{X\inv(g\inv)^E}=\norm{g^E X}$, where the last follows from
  Lemma~\ref{lem:equivalentNorms}.  So, the result follows from
  Theorem~\ref{thm:specificReductionBySupport} by considering $Z\inv$
  instead of~$Z$.
\end{proof}
\end{sloppypar}

\begin{lemma}
\label{lem:expectedNormReduction-2}
Under the assumptions of Lemma~\ref{lem:expectedNormReduction},
and assuming $m\ge\Pr(W=g)$ for all $g\in G$,
\[
\E_{g\in W}(\norm{Xg^E}^2)
        \le \left(\bigl(\Pr(E=0)\bigr)^2
                                + \bigl(\Pr(E=1)\bigr)^2\right) \norm{X}^2 +
          2m \Pr(E=0) \Pr(E=1)
\]
\end{lemma}
\begin{proof}
The lemma follows from Lemma~\ref{lem:expectedNormReduction} and
$\Pr(XW=h)=
  \sum_{g\in G}\Pr(W=g)\*\Pr(X=hg\inv)\le m\sum_{g\in G}\Pr(X=hg\inv)=m$.
\end{proof}

% \begin{theorem}
% \label{thm:reduction}
% Let $X$ and $W$ be random variables on a group~$G$.  Let $E$ be a
% uniform $\{0,1\}$-random variable and let $X$ and $E$ be independent.
% Assume a positive constant~$c$ such that $\norm{X}^2\ge c\max_{g\in
%   G}\Pr(W=g)$.  Then for $f(g)=\norm{Xg^E}$,
% \begin{eqnarray*}
% && \PrOP\left(f(W)
%        < \alpha\sqrt{\lambda}\norm{X}\right)\ge
%    1-\frac{1}{\lambda},\\
% && \quad\mbox{for\ }\alpha =
%         \sqrt{\bigl(\Pr(E=0)\bigr)^2 + \bigl(\Pr(E=1)\bigr)^2
%                     + \frac{2}{c^2} \Pr(E=0) \Pr(E=1)}.
% \end{eqnarray*}
% \end{theorem}
% \begin{proof}
% We apply Theorem~\ref{thm:reduction}.
% Since $E$ is uniform, the value of $\alpha$ reduces to
% $\alpha=\sqrt{(1+1/c)/2}$.
% \end{proof}

The estimate of the next theorem is used for Case~1 in
Section~\ref{sec:algorithm}.

% \begin{theorem}
% \label{thm:reduction}
% Let $X$ and $W$ be random variables on a group~$G$.  Let $E$ be a
% uniform $\{0,1\}$-random variable and let $X$ and $E$ be independent.
% Assume a positive constant~$c$ such that
% $\norm{X}^2\ge c\max_{h\in G}\Pr(W=h)$.
% Let $\lambda>1$ and let $g\in G$ be drawn from the distribution of~$W$.
% Then with
% probability at least $1-1/\lambda$,
% \[
%    \norm{Xg^E} < \sqrt{\frac{\lambda}{2}\left(1+\frac{1}{c}\right)}\,\norm{X}.
% \]
% \end{theorem}
% \begin{proof}
% Let $\alpha>0$ satisfy $\alpha^2 =
%         \bigl(\Pr(E=0)\bigr)^2 + \bigl(\Pr(E=1)\bigr)^2
%                     + \frac{2}{c^2} \Pr(E=0) \Pr(E=1) = 1/2+1/(2c)$.
% Let $m=\max_{g\in G}\Pr(W=g)$ in Lemma~\ref{lem:expectedNormReduction-2}.
% Note that $c\le\norm{X}^2/m$.
% Combining Lemma~\ref{lem:expectedNormReduction-2} with the hypotheses of the
% theorem yields for $f(g)=\norm{Xg^E}$,
% $\E(f(W)) \le \alpha^2 \norm{X}^2$.
% By Markov's inequality (Lemma~\ref{lem:markov}),
% this implies
% $\PrOP\big(f(W)\ge\lambda \alpha^2\norm{X}^2\big)\le
% 1/\lambda$.
% This is equivalent to
% \linebreak[0]$\PrOP\left(f(W)
%        < \alpha\sqrt{\lambda}\norm{X}\right)\ge
%    1-\frac{1}{\lambda}$,
% from which the theorem follows.
% \end{proof}

\begin{theorem}
\label{thm:reductionBoundedSupport}
Let $X$ and $W$ be random variables on a group~$G$.  Let $E$ be a
uniform $\{0,1\}$-random variable and let $X$ and $E$ be independent.
Assume $\Pr(W\notin{\supp(X)})\ge\delta$.  Assume further that
$\Pr(W=g)=\max_{h\in G}\Pr(W=h)$ for all $g\in\supp(X)$.
Let $\lambda>1$ and let $g\in G$ be drawn from the distribution of~$W$.
Then with
probability at least $1-1/\lambda$,
\[
   \norm{Xg^E} < \sqrt{\lambda(1-\delta/2)}\,\norm{X}.
\]
\end{theorem}
\begin{proof}
Let $m=\max_{g\in G}\Pr(W=g)$.  Then $m|\supp(X)|+\delta\le 1$.
So $|\supp(X)|\le(1-\delta)/m$.
Next, $\norm{X}^2\ge1/|\supp(X)|\ge m/(1-\delta)$ by
Lemma~\ref{lem:minNorm-ell2}.  So, $m\le(1-\delta)\norm{X}^2$.
Combining this inequality with Lemma~\ref{lem:expectedNormReduction-2}
and $\Pr(E=0)=1/2$ yields
$\E_{g\in W}(\norm{Xg^E}^2)\le(1/2)\norm{X}^2+m/2
\le (1-\delta/2)\norm{X}^2$.
By Markov's inequality (Lemma~\ref{lem:markov}),
this implies for the function $f(g)=\norm{Xg^E}^2$ from~$G$ to the
real numbers, that
$\PrOP\big(f(W)\ge\lambda (1-\delta/2)\norm{X}^2\big)\le
1/\lambda$.
This is equivalent to
$\PrOP\!\left(f(W)
       < \sqrt{\lambda(1-\delta/2)}\,\norm{X}\right)\ge
   1-\frac{1}{\lambda}$,
from which the theorem follows.
\end{proof}

\IGNORE{
{\bf Choose $c=2$ and $\lambda=4/3$ and then
$\sqrt{\frac{\lambda}{2}\left(1+\frac{1}{c}\right)} = 5/6$
and $1-1/\lambda=1/4$.  Can we make it better?}
} % end IGNORE

% end IGNORE

\section{Decomposition of a random variable}
\label{sec:decomposition}

% {\bf ALTERNATIVE METHOD:  In this section, can we decompose
% $\calR_ig_i^{E_i}$ into the identically distributed random variable
% $(X^JY^{1-J})g_i^{E_i}
% =(Xg_i^{E_i})^J\,(Yg_i^{E_i})^{1-J}$, and then argue with
% $W=(Yg_i^{E_i})^{1-J}$?  Is there a way to make this work?
% }

\IGNORE{
{\Large\bf Consider $W\rightarrow T$ in this section, and in
references to it in Section 6 (where referencing $W=g^{E_i}$).
But this might conflict with $W=W'^KT^{1-K}$.
}
}

One key to this paper is that given random variables~$Z$ and~$X$, we
can decompose $Z$ into~$X$ and a new random variable~$Y$, subject to a
certain ``domination condition''.
In this section, the variable~$Z$ plays the role of $\calR_i$ in the
main algorithm, and the variable~$W$ plays the role of $g_i^{E_i}$ in
the main algorithm.  Hence in the application to the main algorithm,
$W$ can have only two values, $g_i$ and the identity element.
Further, $\Pr(W=g_i)=\Pr(E_i=1)$.

\begin{definition}
For $X$ and $Y$ random variables on a group~$G$, the statement
$X\probeq Y$ means $\forall g\in G,\,\Pr(X=g)=\Pr(Y=g)$
(i.e.~$X$ and~$Y$ are identically distributed).
\end{definition}

\begin{lemma}[decomposition]
\label{lem:decomposition}
Let $Z$ and $X$ be random variables on a group~$G$ and let $I$ be a
$\{0,1\}$-random variable with $I$ independent of~$X$.  Assume that
$\Pr(I=1)\Pr(X=g)\le\Pr(Z=g)$ for all $g\in G$.  Then there is a
decomposition of~$Z$ such that $Z\probeq X^I Y^{1-I}$ for
all~$g\in G$, where $Y$ is a random variable on~$G$
independent of~$I$ and is unique up to probability density.
\end{lemma}
\begin{proof}
Choose $Y$
independent of~$X$ and~$I$ to have a probability density function satisfying
$\Pr(I=0)\Pr(Y=g)=\Pr(Z=g)-\Pr(I=1)\Pr(X=g)$.
\end{proof}

In Section~\ref{sec:algorithm}, this lemma will be used repeatedly for
such decompositions as $\calR_i\probeq X^IY^{1-I}$.  This allows us to draw a
$g\in G$ from the distribution of $\calR_i$ with the knowledge that
with probability $\Pr(I=1)$, it is as if the group element had been
drawn from the distribution of~$X$.  Since we can choose $X$
arbitrarily subject to the domination
condition~$\Pr(I=1)\Pr(X=g)\le\Pr(\calR_i=g)$, this gives us a lot of
flexibility.

Once an $\varepsilon$-uniform random variable is available for some
$\varepsilon<1$, the next lemma shows how to iterate to improve the
uniformity.

\begin{lemma}
\label{lem:acceleratorUniform}
Let $X$ and $Y$ be independent random variables on a group~$G$.  Let
$X$ be $\delta$-uniform and $Y$ be $\varepsilon$-uniform.  Then $XY$
is a $\delta\varepsilon$-uniform random variable.
\end{lemma}
\begin{proof}
Let $U$, $V$, $I$ and $J$ be independent random variables.  Let $U$
and $V$ be uniform on~$G$ and let $I$ and $J$ be on $\{0,1\}$, where
$\Pr(I=0)=\delta$ and $\Pr(J=0)=\varepsilon$.  Further, by
Lemma~\ref{lem:decomposition}, we can write $X\probeq U^IA^{1-I}$ and
$Y\probeq V^JB^{1-J}$ for some random variables, $A$ and~$B$.  
Note that $\Pr(I=0)\Pr(A=g)\le2\delta/|G|$ and
$\Pr(I=1)\Pr(U=g)=(1-\delta)/|G|$ for all $g\in G$ and similarly for
$J$, $V$, $B$ and~$\varepsilon$.
By Lemma~\ref{lem:uniform},
$UV$, $UB$ and $AV$ are all uniform.  So there is a $\{0,1\}$-uniform
random variable~$W$ and a $\{0,1\}$-random variable~$K$ such that
$XY\probeq W^K(AB)^{1-K}$ with $\Pr(K=0)=\Pr(I=0\mbox{ and }J=0)
=\delta\varepsilon$.
So $\Pr(XY=g)\ge\Pr(K=1)\Pr(W=g)=({1-\delta\varepsilon})/|G|$.
Also $\Pr(K=0)\Pr(AB=g)\le2\delta\varepsilon/|G|$.
So $\Pr(XY=g)\le\Pr(K=1)\*\Pr(W=g)+\Pr(K=0)\*\Pr(AB=g)
=({1+\delta\varepsilon})/|G|$.
\end{proof}

\IGNORE{

{\bf DO WE EVER USE THIS NEXT LEMMA?}
\begin{lemma}
\label{lem:decomposition-ell2}
Let $X$, $Y$, $Z$ and $I$ be independent random variables
with $X$, $Y$ and $Z$ on the group~$G$, and $I\in\{0,1\}$.
Then $\norm{XY^IZ^{1-I}}\le\norm{XY}\Pr(I=1)+\norm{XZ}\Pr(I=0)$.
\end{lemma}
\begin{proof}
Then $\norm{XY^IZ^{1-I}}=\norm{(XY)^I(XZ)^{1-I}}
=\norm{\varphi(XY)\Pr(I=1)+\varphi(XZ)\Pr(I=0)}_2
\le\norm{XY}\Pr(I=1)+\norm{XZ}\Pr(I=0)$.
\end{proof}

} % end IGNORE

The next lemma from linear algebra is a standard calculation on
vectors in the $\ell^2$ norm.  It is needed to prove the succeeding
Theorem~\ref{thm:nextStep}.  The $\ell^2$ vectors of the lemma will
correspond to vectors of dimension~$|G|$, where a $G$-valued random
variable is considered as $(\Pr(X=g_1),\ldots,\Pr(X=g_{|G|}))$.

\begin{lemma}
\label{lem:norm-comparison-ell2}
Let $c$ and $\alpha$ be constants.  Let $x$, $y$, $x'$ and~$y'$ be
vectors in the $\ell^2$ norm.
Assume $\norm{x'}_2\le\alpha\norm{x}_2$, $\norm{y'}_2\le\norm{y}_2$, and
$\norm{x}_2\ge c\norm{y}_2$
for $0<\alpha\le 1$ and $c>0$.
Then $\norm{x'+y'}_2
\le \left(({1+\alpha c})/\sqrt{1+c^2}\right)\,\norm{x+y}_2$.
\end{lemma}
\begin{proof}
Note that $\norm{x}_2+\norm{y}_2\le
\left((1+c)/\sqrt{1+c^2}\right)\sqrt{\norm{x}_2^2+\norm{y}_2^2}$.
Let $d=(c-\alpha c)/(1+c)$.
The proof follows from $\norm{x'+y'}_2\le\norm{x'}_2+\norm{y'}_2
\le\alpha\norm{x}_2+\norm{y}_2
\le(\alpha+d/c)\norm{x}_2+(1-d)\norm{y}_2
=\left((1+\alpha c)/(1+c)\right) \*\, (\norm{x}_2+\norm{y}_2)$.
Note that for fixed $\norm{x}_2$ and $\norm{y}_2$, $\norm{x+y}_2$ is
minimized when $x$ and $y$ are perpendicular.
In this case, define $c'$ such that $\norm{x}_2=c'\norm{y}_2$
and observe that
$\left((1+\alpha c)/(1+c)\right) \*\, (\norm{x}_2+\norm{y}_2)
=\left(({1+\alpha c})/(1+c)\right) \*
        ((c'+1)/\sqrt{{c'}^2+1}) \*
                                \sqrt{\norm{x}_2^2+\norm{y}_2^2}
\le\left(({1+\alpha c})/\sqrt{1+c^2}\right) \*
                                \sqrt{\norm{x}_2^2+\norm{y}_2^2}
=\left(({1+\alpha c})/\sqrt{1+c^2}\right)\*\,\norm{x+y}_2$.
\end{proof}

The estimate of the next theorem is used for Cases~2 and~3 in
Section~\ref{sec:algorithm}.

\IGNORE{
{\bf ALTERNATIVE:  $\norm{ZW}\le \alpha\norm{X}+\norm{Y}
\le \beta(\norm{X}-\norm{Y})\le\beta\norm{Z}$,
where
$\beta\ge(\alpha\norm{X}+\norm{Y})/(\norm{X}-\norm{Y})
= (\alpha c+1)/(c-1)$.  These estimates assume first Y in same
direction as X, and then Y in opposite direction from X.  Stuff below
yields finer estimate.  Can we do it more compactly?  In particular,
it's impossible that Y be opposite to X, since that would make Y
negative.
The largest angle implies that Y is perpendicular to X.}
} % end IGNORE

\begin{theorem}
\label{thm:nextStep}
Let $X$, $Y$, $Z$ and~$W$ be random variables on a group~$G$ and let
$J$ be a $\{0,1\}$ random variable.  Let $X$, $Y$, $W$ and $J$ be independent,
and let $Z\probeq X^JY^{1-J}$.  Let $\norm{XW}\le\alpha\norm{X}$
and $\Pr(J=1)\norm{X}\ge c\Pr(J=0)\norm{Y}$ for some $0<\alpha\le 1$ and $c>0$.
Then
\[
\norm{ZW}\le\left(({1+\alpha c})/\sqrt{1+c^2}\right)\,\norm{Z}.
\]
\end{theorem}
\begin{sloppypar}
\begin{proof}
By Lemma~\ref{lem:bounded-ell2}, $\norm{YW}\le\norm{Y}$.
By Lemma~\ref{lem:norm-comparison-ell2},
$\norm{ZW}=\norm{\Pr(J=1)\varphi(XW)+\Pr(J=0)\varphi(YW)}_2
\le \left(({1+\alpha c})/\sqrt{1+c^2}\right)\*
        \norm{\Pr(J=1)\varphi(X)+\Pr(J=0)\varphi(Y)}_2
= \left(({1+\alpha c})/\sqrt{1+c^2}\right) \norm{Z}$.
\end{proof}
\end{sloppypar}

\IGNORE{
% NOT NEEDED:
\begin{remark}
\label{rem:nextStepConstants}
  If $c>2\alpha/(1-\alpha^2)$, then $\norm{ZW}<\norm{Z}$ in
  Theorem~\ref{thm:nextStep}.  For $\alpha$ close to~1, and
  $c=k/(1-\alpha)$, $\left(({1+\alpha c})/\sqrt{1+c^2}\right) =
  1-(1-1/k)(1-\alpha)+O((1-\alpha)^2)$.
%  PROOF: Let $\epsilon=1-\alpha$ and take power series in $\alpha$.
\end{remark}
}

\IGNORE{
{\bf Remove previous lemma, theorem and remark if we go with next theorem?}
} % end IGNORE

The estimate of the next theorem is used for Case~1 in
Section~\ref{sec:algorithm}.

\begin{theorem}
\label{thm:newNextStep}
  Let $X$, $Y$, $Z$ and~$W$ be random variables on a group~$G$ and let
  $J$ be a $\{0,1\}$ random variable.  Let $X$, $Y$, $W$ and $J$ be
  independent, and let $Z\probeq X^JY^{1-J}$.  Assume constants $m$,
  $c$ and~$\alpha$ satisfying the following.  Let
  $\Pr(J=0)\Pr(Y=g)\le m$ for all $g\in G$, and further let
  $\Pr(J=0)\Pr(Y=g)=m$ when $g\in\supp(X)$.  Let
  $\norm{XW}\le\alpha\norm{X}$ for some $0<\alpha\le 1$.
  Let $c= (1-\alpha^2)\,(\Pr(J=1))^2 /
    ( m|A|\Pr(Z\notin A)+1 )$ for $A=\supp(X)$.
Then
\[
\norm{ZW}\le\sqrt{1-c}\,\norm{Z}
\]
\end{theorem}
\begin{sloppypar}
\begin{proof}
Note that $\Pr(J=0)\Pr(YW=g)=\Pr(J=0)\sum_{h\in G}\Pr(Y=h)\Pr(W=h\inv g)
\le$\linebreak[0]$m\sum_{h\in G}\Pr(W=h\inv g)=m$.
Since $\Pr(X=g)=0$ for $g\notin\supp(X)$ and $\Pr(J=0)\Pr(Y=g)=m$
for $g\in\supp(X)$,
we have
\begin{eqnarray*}
\lefteqn{\sum_{g\in G}2\Pr(J=1)\Pr(J=0)\Pr(XW=g)\Pr(YW=g)}\qquad\qquad \\
 &\le&2m\Pr(J=1)\sum_{g\in G}\Pr(XW=g) \\
 &=&2m\Pr(J=1) \\
 &=&2m\Pr(J=1)\sum_{g\in \supp(X)}\Pr(X=g) \\
 &=&2\Pr(J=1)\Pr(J=0)\sum_{g\in G}\Pr(X=g)\Pr(Y=g).
\end{eqnarray*}
By Lemma~\ref{lem:bounded-ell2}, $\norm{YW}\le\norm{Y}$.
Hence,
\begin{eqnarray*}
\norm{ZW}^2&=&\norm{\Pr(J=1)\varphi(XW)+\Pr(J=0)\varphi(YW)}_2^2 \\
&=&\sum_{g\in G}(\Pr(J=1)\Pr(XW=g)+\Pr(J=0)\Pr(YW=g))^2 \\
&=&(\Pr(J=1)\norm{XW})^2+(\Pr(J=0)\norm{YW})^2 + \\
 && \qquad    \sum_{g\in G}2\Pr(J=1)\Pr(J=0)\Pr(XW=g)\Pr(YW=g) \\
&\le&\alpha^2(\Pr(J=1)\norm{X})^2+(\Pr(J=0)\norm{Y})^2 + \\
 && \qquad    \sum_{g\in G}2\Pr(J=1)\Pr(J=0)\Pr(X=g)\Pr(Y=g) \\
&=&\norm{Z}^2 - (1-\alpha^2)\*\,(\Pr(J=1)\norm{X})^2 \\
&\le&(1-c)\norm{Z}^2
\end{eqnarray*}
providing $({1-\alpha^2})\*\,\left(\Pr(J=1)\norm{X}\right)^2\ge c\norm{Z}^2$
for $c>0$.

We find such a~$c$.  Let $A=\supp(X)$.  Note that
$m|A|+\Pr(J=1)=\Pr(J=0)\Pr(Y\in A)+\Pr(J=1)\le1$.  One can show that
$\norm{X}^2/\norm{Z}^2$ is minimized when $X$ is uniform on~$A$.  In
this case, $\norm{X}^2=1/|A|$ and we have
$\norm{Z}^2\le m^2(\Pr(Z\notin A)/m)+|A|(m+\Pr(J=1)/|A|)^2
\le m\Pr(Z\notin A)+1/|A|$.  So, if
$c= ({1-\alpha^2})\*\,(\Pr(J=1))^2 /\mathbreak ( m|A|\Pr(Z\notin A)+1 )$, then
$({1-\alpha^2})\*\,\left(\Pr(J=1)\norm{X}\right)^2\ge c\norm{Z}^2$.
\end{proof}
\end{sloppypar}

\section{Fuzzy subgroups and escaping from a set}
\label{sec:escape}

Section~\ref{sec:reduction} constructs an element $g\in G$ such that
$\norm{\calR_i g^{E_i}}<c\norm{\calR_i}$ for some $c<1$.  However, it
fails if, for example, $X$, $W$ and $XW$ all have identical
distribution.  This is part of a larger class of examples.  If $\calR_i$
is the uniform distribution on a proper subgroup $H<G$, then any
construction of~$X$ and~$W$ from~$\calR_i$ will fail to produce a $g\in
G$ with $\norm{\calR_i g^{E_i}}<c\norm{\calR_i}$, since $\norm{\calR_i}$
is already minimized among random variables on $H<G$.
Hence, when the methods of Section~\ref{sec:reduction} fail, we must
demonstrate that this implies that $\calR_i$ is close to a uniform
distribution on a proper subgroup $H<G$.

The following surprising lemma is the key.  It shows that if one
cannot escape a set~$A=B\inv B$ with reasonable probability simply by
multiplying two random elements of the set~$A$, then one must be
``stuck'' in a proper subgroup.  Loosely speaking, either the product
of two random elements of~$A$ ``escapes'' from the set~$A$, or else
$A$ must be a ``fuzzy subgroup'' of~$G$ in the sense that $A$ is close
in probability to some subgroup of~$G$.  In the latter case, we use
the generators of~$G$ to construct a $g\notin H$ so that $\calR_ig^{E_i}$
``escapes'' the set~$H$.

The proof proceeds by constructing a multiplication table for products
of elements of~$A$.  If $gh\notin A$ for $g,h\in A$, then we think of
$gh$ as a ``hole'' in the multiplication table.  We then augment the
multiplication table to include $gh$, and show that the number of
holes in the multiplication table for $A\cup{gh}$ has been reduced.

\begin{lemma}
\label{lem:key}
Let $A\subseteq G$ satisfy $A=A\inv$.  Let $\delta<1/4$.  Assume
$\forall g\in A, |Ag\setminus A|\le\delta|A|$.  Then
\[ |AA\setminus A|\le{\delta\over 1-2\delta}|A|.
\]
Furthermore, $AA$ is a subgroup of~$G$.
\end{lemma}
\begin{proof}
Define $\phi(g)=|\{a\in A\colon~ag\notin A\}|$.  The hypothesis can
be re-phrased as $\forall g\in A, \phi(g)\le\delta|A|$.  From this it
follows that
\[\forall g,h\in A, \phi(gh)\le\phi(g)+\phi(h)\le 2\delta|A|.
\]
Similarly, for all $a,b,c,d\in A$,
$\phi(abc)\le3\delta|A|$ and $\phi(abcd)\le4\delta|A|$.

So, for $g,h\in A$ such that $gh\notin A$, there are at least
$(1-2\delta)|A|$ pairs $(u,v)$ such that $gh=uv$.  To see this, note
that $v=u\inv gh$ and so we are counting the number of pairs $(u,v)\in
A\times A$ such that $u\inv gh\in A$.  This number is
$|A|-\phi(gh)\ge(1-2\delta)|A|$.
Hence, $|AA\setminus A|\le \delta/(1-2\delta)|A|$ as required by the
lemma.

It remains to show that $AA$ is a group.  Since it is closed under
inverses, we must show that it is closed under multiplication.  For
$a,b\in A$, it is clear that $ab\in AA$.  Given $a\in A$ and $c,d\in
AA\setminus A$, we must demonstrate membership in three cases:
$ac,ca,cd\in AA$.

We first show that $cd\in AA$, where $c=gh$, $d=uv$, $c,d\in
AA\setminus A$ and $g,h,u,v\in A$.  Note that
$\phi(ghuv)\le4\delta|A|<|A|$.  Therefore,
$\exists w\in A$ such that $wghuv\defeq x\in A$.  So,
$cd=(gh)\,(uv)=w\inv x\in AA$.

A similar argument holds to show that $ca\in AA$, where $c=gh$, $c\in
AA\setminus A$ and $c,g,h\in A$.  It follows from noting that
$\phi(gha)\le3\delta|A|<|A|$.  Finally, $ac=c^{-1}a^{-1}$, and so the
case of $ac$ reduces to the previously proved case of $ca$.
\end{proof}

\begin{remark}
Examination of the proof shows that the hypothesis could be weakened
to $\delta<1/3$ or further, at the cost of showing that $A^k$ is a
group for some sufficiently large~$k$.
\end{remark}

\IGNORE{
{\bf
Page 11, remark 2. In fact if delta < 1/3 then you can prove
that A^2 is a group, don't need A^k.
}
}

One interpretation of Lemma~\ref{lem:key} is that for random $(u,v)\in
A\times A$, $uv\notin A$ with some constant probability, or $AA$ is
close to a group and so $ug\notin A$ with some constant probability
for some group generator~$g$.

\IGNORE{

{\bf Probabilistic version:  Let $A$ be a multiset (i.e. an unordered
list in which elements of $G$ may be repeated).  The above proof goes
through as before.  CHECK THIS!!  Now some count in which we have
$\delta|A|$ gets reinterpreted as
a probability, $\delta$.

Maybe want to note hypothesis of Lemma~\ref{lem:key} can be
interpreted as $\Pr(X_A g\notin A)\le\delta \forall g\in A$ where
$X_A$ is a uniform random variable on~$A$, and $\Pr(X_A\notin A)=0$.
This generalizes to hypotheses of $\Pr(Y\in A)=1$ and $\Pr(Yg\notin
A)\le\delta\forall g\in $, where $Y$ has an arbitrary distribution.
The conclusion is then $AA$ is a group and $\Pr(Y\times Y\notin
B)\le\delta$ where $B=\{(u,v)\colon~uv\in A\}$.
}

}

\begin{theorem}
\label{thm:escape}
Let $k>1$ and $0<\varepsilon<1$ be arbitrary constants and let
$\delta=(2+k^2\varepsilon)/(k-2)$.  Assume $\delta\le1/4$.
Let $A\subseteq G=\gen{S}$ satisfy
$A=A\inv$.  Then one of the following is true.
\begin{enumerate}
\item Given a random $(u,v)\in A\times A$ drawn from a uniform
distribution, $uv\notin A$ with probability at least $\varepsilon$.
\item $\exists A'\subseteq A$ with $|A\setminus A'|<2|A|/k$ such
that $A'A'$ is a subgroup of~$G$.
Furthermore
\[
 |A'A'\setminus A'|\le{\delta\over 1-2\delta}|A'|.
\]
\end{enumerate}
\end{theorem}
\begin{proof}
Define $B\subseteq A\times A$ such that $B=\{(g,h)\colon~gh\in A\}$.
If $|A\times A\setminus B|>\varepsilon|A\times A|$, then a random $(u,v)\in
A\times A$ satisfies $uv\notin A$ with probability at least
$\varepsilon$, and we are done.

Otherwise, $|A\times A\setminus B|\le\varepsilon|A\times A|$.  Note
that $k\varepsilon<1$.
Let
\[ A'=\{g\in A\colon~|Ag\setminus A|\le
k\varepsilon|A|\mbox{\rm\ and\ }|gA\setminus A|\le k\varepsilon|A|\}.
\]
Note that $A'={A'}\inv$.
Also $|A\setminus A'|<2|A|/k$.
To see the latter, note that $|\{g\colon~|Ag\setminus A|>k\varepsilon|A|\}|
 < |\{(u,g)\colon~ug\notin A\}|/(k\varepsilon|A|)
 = |A\times A\setminus B|/(k\varepsilon|A|) \le |A|/k$.

% {\bf REVISED PARAGRAPH: }

Therefore $|(A'g\cap A)\setminus A'|\le|A\setminus A'|\le(2/k)|A|$.
Also $|A'g\setminus A|=|g\inv A'\setminus A|\le k\varepsilon|A|$
       for all~$g\in A$.
Hence $|A'g\setminus A'|\le
  (2/k)|A|+k\varepsilon|A|
  \le ((2+k^2\varepsilon)/k)|A|$.
But $|A|<|A'|/(1-2/k)$ follows from $|A|-|A'|=|A\setminus A'|<2|A|/k$.
(The coefficient $1/(1-2/k)$
  is positive since $\delta\le 1/4$ implies that $k\ge 10$.)
Hence $|A'g\setminus A'| < ((2+k^2\varepsilon)/(k-2))\,|A'| = \delta|A'|$
  for all $g\in A'$.

Since $\delta\le 1/4$ and
$|A'g\setminus A'|\le\delta|A'|$ for all~$g\in A'$, we can invoke
Lemma 5.1 on $A'$ and conclude that $A'A'$ is a subgroup
of~$G$.  The bounds on $|A'A'\setminus A'|$ follow from the same lemma.
\end{proof}

% \begin{quotation}
% \noindent
% {\em
% Explanation of {\rm ``Then $|A'|\ge ((k-1)/k)|A|$.''
% or now: ``Then $|A\setminus A'|<2|A|/k$''}
% \newline
% $|A\times A\setminus B|=\mbox{\# holes}$ and $|A'|=\#$ rows with
% $k\varepsilon|A|$ holes or less (a row for each~$g$).
% $|A\setminus
% A'|\le 2{\mbox{\rm\# holes} \over k\varepsilon|A|}
% =2{|A\times A\setminus B|\over k\varepsilon|A|}
% \le2{\varepsilon|A\times A|\over k\varepsilon|A|}\le{2|A|\over k}$.
% \newline
% This implies ${|A'|\over|A|}\ge1-{2\over k}={k-2\over k}$, which
% implies ${|A\setminus A'|\over|A|}\le{2\over k}$ or $|A'|\ge ((k-2)/k)|A|$.
% }
% \end{quotation}

\IGNORE{

{\bf
  Refinement to try later:  If many holes in $(A\setminus A')\times
  (A\setminus A')$ instead of $A\times (A\setminus A')$, then get finer
  estimate, since we can ignore holes in $(A\setminus A')\times
  (A\setminus A')$.
}

} % end IGNORE

%DELETE THIS COMMENT WHEN PROOF IS STABLE.
% We have $A'A'$ is a group, and $|A'A'|/|A|\ge|A'|/|A|=(k-1)/k=0.9$.
% Therefore, for random $u\in A$, $u\in A'$ with probability at least
% 0.9.  With probability at least 0.5, $|A'A'r\setminus
% A'A'|\ge|A'A'|/2$.  {\bf WANT $ur\notin A$ or else $uv\notin A$ for
% $(u,v)\in A'\times A'$.  $A'A'$ can't move far outside of~$A$ by
% construction.  So, if $A'A'r$ moves far outside of $A$, then there
% must be many $u\in A'$ that are being moved outside of~$A$.}

% $\varepsilon=(\delta(k-2)-2)/k^2$
% $d\varepsilon/dk=\delta/k^2-(2(\delta(k-2)-2))/k^3=0$
% $k\delta+4-2\delta(k-2)=0$
% $4+4\delta-\delta k=0$
% $k=4/\delta+4$

%$\delta<1/4 => (1+x(x-1)\varepsilon)/(x-1)<1/4$, where $|Ag\setminus A|\le
%  x\varepsilon|A|$.  This implies $4x^2 \varepsilon -4\varepsilon -x +5\le0$

\begin{corollary}
\label{cor:escape}
Assume $A$, $k$ and $\varepsilon$ as in Theorem~\ref{thm:escape}.
Let $p=(1/2)-(1/k)$ and
let $r$ be a random subproduct on~$S$.  Assume $(u,v)\in
A\times A$ drawn from a uniform distribution.
Let $g=uv^Ir^{1-I}$ for $I$ a
$\{0,1\}$ random variable with $\Pr(I=1)=p/(p+\varepsilon)$.
Then $g\notin A$ with probability at least
$p\varepsilon/(p+\varepsilon)
> \varepsilon-2k\varepsilon^2/(k-2)$.
\end{corollary}
\begin{proof}
Theorem~\ref{thm:escape} tells us that $uv\notin A$ with probability
at least~$\varepsilon$ or $A'A'$ is a subgroup of~$G$ with
$|A\setminus A'|<2|A|/k$.  In the latter case, $r\notin A'A'$ with
probability at least $1/2$.  Hence, with probability at least 1/2, for
$h\in A'$, $hr\notin A'A'\supseteq A'$.  For $u$ drawn at random from
$A$, $ur\notin A$ with probability at least
$(1/2)|A'|/|A|=(1/2)-(1/k)=p$.

Let $g=uv^Ir^{1-I}$.
Then $\Pr(g\notin A) \ge \min(\varepsilon\Pr(I=1),p\Pr(I=0))
 =p\varepsilon/(p+\varepsilon)
 = \varepsilon(k-2)/(2k\varepsilon+k-2)
 > \varepsilon-2k\varepsilon^2/(k-2)$.
\end{proof}

\begin{remark}
\label{rem:constants}
  Consider the equation $\delta=(2+k^2\varepsilon)/k$ of
  Theorem~\ref{thm:escape}.  The variable $\varepsilon$ is maximized
  when $k=4/\delta+4$.  Taking $\delta=1/4$ implies $k=20$ and
  $\varepsilon=1/160$ when it is maximized.  In this case,
  Corollary~\ref{cor:escape} produces a $g\notin A$ with probability
  at least $\varepsilon>0.006$.
\end{remark}

% DELETE THIS WHEN PROOF IS STABLE.
% {\bf
% This is approximate.
% $|A'A'|\ge\left({k-1\over k}\right)^2|AA|$.\newline
% $|A'A'r\setminus A'A'|\ge|A'A'|/2\ge\left({k-1\over
%     k}\right)^2|A|^2/2$.\newline
% $|A'A'|\ge\left({k-1\over k}\right)^2|A|^2$
% \newline
% overflow $\Rightarrow |A'A'r\setminus AA|\ge|AA|$.
% }

\section{Fibonacci Cube algorithm for semi-uniform random generation}
\label{sec:algorithm}

We now have all of the algorithmic components outlined in
Section~\ref{sec:outline}.  The goal of this section is only to
construct $g_i$ for which $g_1^{E_1}\cdots g_t^{E_t}$ is
{\em semi-uniform}.

\begin{definition}
A random variable on a group~$G$ is an $\varepsilon$-semi-uniform random
variable if $\Pr(X=g)\ge{1/|G|-\varepsilon/|G|}$ for all $g\in G$.
The random variable is semi-uniform if it is
$\varepsilon$-semi-uniform for some~$\varepsilon>0$.
\end{definition}

\subsection{Algorithm}
\label{sec:pseudo-code}
Given a random variable~$\calR_i$ on a
group~$G=\gen{S}$, we wish to construct $g_i\in G$ such that
$\norm{\calR_i g_i^{E_i}}/\norm{\calR_i}<c<1$ for some constant~$c$
and for $E_1,E_2,\ldots$ independent uniform $\{0,1\}$-random
variables.  By Lemma~\ref{lem:bounded-ell2},
$\norm{\calR_ih^{E_i}}\le\norm{\calR_i}$ for all $h\in G$.  Hence, we
will construct $g_i$ that has only some constant probability of
satisfying $\norm{\calR_i g_i^{E_i}}/\norm{\calR_i}<c<1$.  We then set
$\calR_{i+1}=\calR_i g_i^{E_i}$, knowing that
$\norm{\calR_{i+1}}=\norm{\calR_i g_i^{E_i}}$, even if $g_i$ did not
succeed.  We can then try again by constructing $g_{i+1}$.  In Case~2
below, we define $\calR_{i+1}=g_i^{E_i}\calR_i$ instead of
$\calR_{i+1}=\calR_i g_i^{E_i}$, but this does not change the spirit
of the algorithm.

We call the algorithm below the Fibonacci Cube algorithm by allusion to the
Fibonacci series.  Like the Fibonacci series, each group element is
derived from the previous elements of the series.  It is a cube
algorithm since $\calR_i=h_1^{E_1}\cdots h_k^{E_k}$ for exponents that
are independent uniform $\{0,1\}$-random variables.
The pseudo-code for the algorithm is simple.
\\
{ % This must be outermost level in TeX.
\obeylines \obeyspaces
{\obeyspaces\global\let =\ }% let active space = control space
% {\obeyspaces\gdef {\phantom{X}}}% let active space = control space
{\bf{Algorithm Fibonacci-Cube}}
INPUT:  Black box group $G=\gen{S}$
OUTPUT: $\calR_t\inv \overline\calR_t$ %
for $\overline\calR_t$ an independent copy of $\calR_t$;
                [ For large enough~$t$, %
$\Pr(\calR_t\inv\overline\calR_t=g)\ge(3/4)(1-\beta)^2/|G|$ for all $g\in G$ ]
PARAMETERS:  positive constants $a$, $b$ and~$c$; % required '%'
$\alpha$ and $\rho$ dependent on $a$, $b$ and~$c$
                        such that %
$\norm{\calR_{i+1}}\le \alpha\norm{\calR_i}$ %
with probability at least $\rho$
                        unless $\calR_i\inv\overline\calR_i$ already %
satisfies the conditions on $\calR_t\inv\overline\calR_t$
   Let $\calR_1$ be the identity element with probability~1
   Let $t=\log|G|/\log \alpha^{-2\rho}$
   For $i=1$ to $t-1$
     Let $d=1/a+1/b+1/c$
     Let $j\in 1$, $2$ or $3$ with
       probability $1/(ad)$, $1/(bd)$ or $1/(cd)$, respectively
     Goto Case j
     Case 1:
       Choose $g_i$ from distribution of $\calR$
       Set $\calR_{i+1}=\calR_ig_i^{E_i}$
     Case 2:
       Choose $g_i$ from distribution of $\calR$
       Set $\calR_{i+1}=g_i^{E_i}\calR_i$
     Case 3:
       Choose $g_i$ from distribution of random subproducts on~$S$
       Set $\calR_{i+1}=\calR_ig_i^{E_i}$
   Return $\calR_t\inv \overline\calR_t$ %
for $\overline\calR_t$ an independent copy of $\calR_t$
}
\ \newline

Note that the output of the algorithm is in terms of a  random
variable $\calR_t=h_1^{E_1}\cdots h_t^{E_t}$, where
$(h_1,\ldots,h_t)$ is a reordering of $(g_1,\ldots,g_t)$.  So, an
implementation of the algorithm would need only to record the elements
$(h_1,\ldots,h_t)$.  An element from the distribution
$\calR_t\overline\calR_t$ is then computed as
$(h_t\inv)^{\overline E_t}\cdots(h_1\inv)^{\overline E_1}
(h_1)^{E_1}\cdots (h_t)^{E_t}$ where each of
$E_1,\ldots,E_t,\overline E_1,\ldots,\overline E_t$ is independently equal
to zero or~one with probability~1/2.

The random variable produced by the Fibonacci cube algorithm is used
to produce
a $\gamma$-uniform random element.  One can then use
Lemma~\ref{lem:acceleratorUniform} to produce $\varepsilon$-uniform
random elements for arbitrarily small~$\varepsilon$.

\subsection{Overview of proof}
\label{sec:overview-proof}

The immediate goal is to prove Lemma~\ref{lem:main}, that
$\norm{\calR_{i+1}}\le c\norm{\calR_i}$ with probability at least
$\rho>0$ for some positive $c<1$.

In Cases~1 and~2, $\calR_{i+1}=\calR_i g_i^{E_i}$ or
$\calR_{i+1}=g_i^{E_i} \calR_i$, for $g$ drawn from~$W=\calR_i$.
In Case~3, $\calR_{i+1}=\calR_i g_i^{E_i}$ for $g$ a random
subproduct.  The proof proceeds by decomposing both $\calR_i$ and~$W$
as follows
into products of random variables that are easier to analyze.
\begin{eqnarray*}
  && \calR_i\probeq X^JY^{1-J} \\
  && W=\calR_i\probeq W'^K T^{1-K} \mbox{ (Cases~1 and~2 only)}\\
  && W=W' \mbox{ is a random subproduct on the group generators (Case~3)}
\end{eqnarray*}

The general approach in each case is to define $X$, $J$, $W'$ and $K$
so that $\norm{Xg^{E_i}}\le \a\norm{X}$ for some positive
$\a<1$ and for $g$ drawn from the distribution of~$W'$, with
probability $\rho>0$.  The results of Sections~\ref{sec:reduction}
or~\ref{sec:escape} are used here
(Theorem~\ref{thm:reductionBoundedSupport} for
Case~1, Corollary~\ref{cor:altReductionBySupport} for Case~2, and
Theorem~\ref{thm:escape} for Case~3).

Then a result from Section~\ref{sec:decomposition}
(Theorem~\ref{thm:newNextStep} or Theorem~\ref{thm:nextStep}) is used
to show that $\norm{Xg^{E_i}}\le \a\norm{X}$ implies
$\norm{\calR_i g^{E_i}}\le \beta\norm{\calR_i}$ for some positive
$\beta<1$ and for $g$ drawn from the distribution of~$W'$, with
probability $\rho>0$.

Of course, one wishes to draw $g$ from the distribution of~$W$, rather
than from the distribution of~$W'$.  Since $W=W'^KKT^{1-K}$, a group
element~$g$ drawn from the distribution of~$W$ can be considered to
have been drawn from the distribution of~$W'$ with
probability~$\Pr(K=1)$.  Hence, one observes that the previous result
implies that $\norm{\calR_i g^{E_i}}\le \beta\norm{\calR_i}$ for some
positive $\beta<1$ and for $g$ drawn from the distribution of~$W$,
with probability $\rho\Pr(K=1)>0$.

At any step of the algorithm, one does not know which of the three
cases are satisfied by the current~$\calR_i$.  However, this is not a
problem.  One chooses
the recipe of one of the three cases at random in deciding how to
construct $g_i$ and $\calR_{i+1}$.  If an incorrect case is chosen,
Theorem~\ref{lem:bounded-ell2} guarantees that
$\norm{\calR_{i+1}}=\norm{\calR_i g_i^{E_i}}\le\norm{\calR_i}$.  So,
as long as a correct case is chosen with at least some positive
probability, the algorithm makes progress.

The pseudo-code allows one
to choose positive parameters $a$, $b$ and~$c$ to determine the ratio of
the probabilities for choosing each of the three cases.  However, the
algorithm succeeds with the same asymptotic estimates regardless of
the choice of $a$, $b$ and~$c$.

\subsection{Proof}
\label{sec:proof}

The analysis of the pseudo-code will be in terms of four parameters,
$\beta$, $\delta$ and~$\lambda$, such that $1>\beta>2\delta>0$
and $\lambda>1$.  The parameter values will be chosen based on the
requirements of the proof.

% DELETE when sure don't need this:
%  After having fixed those values, we will
% also require that $|G|>\max(1/\delta,1/(\beta-\delta))$ or,
% equivalently, $\delta>1/|G|$ and $\beta-\delta>1/|G|$.  Hence, we
% assume $|G|$ satisfies $1>\beta>\delta>1/|G|$, $\beta-\delta>1/|G|$,
% and $\lambda>1$.

The analysis of Cases~1, 2 and~3 of this section applies for
$|G|>\max(1/\delta,1/(\beta-\delta))$.  The
analysis finds asymptotic bounds on the time to produce an
$\varepsilon$-uniform random variable on~$G$.  For groups with
order~$|G|\le\max(1/\delta,1/(\beta-\delta))$, one can easily show
that the pseudo-code succeeds in some constant time.

\begin{definition}
Define
\begin{eqnarray*}
 \overline A_x &\defeq&\{g\in G\colon~\Pr(\calR_i=g) > x\}\\
% M   &\defeq&\max_{g\in G}\Pr(\calR_i=g)\\
 m   &\defeq&\min_x\{x\colon~\Pr(\calR_i\notin \overline A_x) > \delta\}\\
\end{eqnarray*}
Note that $m$ and $\overline A_m$ implicitly depend on~$\calR_i$, and
hence on~$i$.  Define $A_m\supseteq \overline A_m$ so that
\[
\forall B\supset A_m,\,
  \Pr(\calR_i\notin B)<\delta\le\Pr(\calR_i\notin A_m).
\]
This need not uniquely define
$A_m$, but any instance satisfying the defining
conditions will suffice.  The condition implies that $A_m$ is maximal
in the sense that $\Pr(\calR_i\notin B)<\delta$ for all $B\supset A_m$.
\end{definition}

\begin{lemma}
\label{lem:AmBounds}
Assume $\max_{g\in G}\Pr(\calR_i=g)\le1-\delta$.
The set $A_m\subseteq G$ satisfies
\[
 \delta \le \Pr(\calR_i\notin A_m) < \delta+m.
\]
Also,
\[
\begin{array}{ll}
  \Pr(\calR_i=g) \ge m & \mbox{\rm for } g\in A_m \\
  \Pr(\calR_i=g) \le m & \mbox{\rm for } g\notin A_m
\end{array}
\]
Further, if $m<\delta$, then
\[
  \delta \le \Pr(\calR_i\notin A_m) < 2\delta
\]
\end{lemma}
\begin{proof}
  The first inequality follows easily from the definition of~$A_m$
  and $\max_{g\in G}\Pr(\calR_i=g)\le1-\delta$.
  For the next two inequalities, note that the definition of
  $\overline A_m$ implies there is a $g\in G$ such that
  $\Pr(\calR_i=g)=m$.  If there were only one such~$g$, one would have
  $A_m=\overline A_m$.  If there are multiple such~$g$, then
  $\Pr(\calR_i=g)=m$ for all $g\in\overline A_m\setminus A_m$.
  The last inequality follows from the first one and $m<\delta$.
\end{proof}

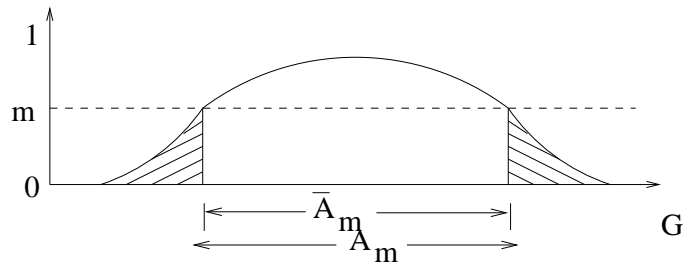
\begin{figure}[htb]
% 80\% reduction from Am.fig
\centerline{\input{Am.pstex_t}}
\caption{Probability density function for $\calR_i$,
         shaded part (outside $\overline A_m$) has area $>\delta$
        and outside of $A_m$, the area is $\ge\delta$ }
\label{fig:Am}
\end{figure}

In the rest of this section, we will isolate a ``Case~0'' to consider
$m\ge\delta$ or $\max_{g\in G}\Pr(\calR_i=g)\le1-\delta$.  In all
other cases, Lemma~\ref{lem:AmBounds} applies with its conclusion that
\[
  \delta \le \Pr(\calR_i\notin A_m) < 2\delta.
\]

\IGNORE{
For intuition, think of $\overline A_x$ as the interior and frontier of an
isocline for a map $g\mapsto\Pr(\calR_i=g)$.
} % end IGNORE

\begin{definition}
Define the random variable $U_{B}$ on~$G$ for a set $B\subseteq G$ by
\[
 \Pr(U_{B}=g) = \left\{
\begin{array}{ll}
 1/|B| &\mbox{\rm for $g\in B$}\\
 0     &\mbox{\rm for $g\notin B$}
\end{array}
\right.
\]
\end{definition}

Recall that $1>\beta>2\delta>0$ and $\lambda>1$ below.
The parameters $\beta$, $\delta$ and
$\lambda$ are fixed throughout.  The parameter~$m$ and the set~$A_m$
depend on~$\calR_i$ and hence on~$i$.
Intuitively, one may think of $1-\beta$ as a
constant against which $m|A_m|$ is measured.  Similarly, one may think
of $\delta$ as a constant against which $\Pr(\calR_i\notin A_m)$ is
measured.  One
thinks of $\Pr(\calR_i\in A_m)-m|A_m|$ as ``large'' if it is larger
than $\beta-\delta$.
In each of the three cases, we will construct $g_i\in
G$ and conclude that there is a $c'<1$ and $\rho'>0$ such that
$\norm{\calR_ig_i^{E_i}}\le c'\norm{\calR_i}$ with probability at
least $\rho'$.

% Intuitively, one may think of $\delta$ as approximating the area under
% the curve outside of $A_m$, on the tails.  One may think of think of
% $m|A_m|$ as a rectangle, and $\beta$ is the ``area'' under $Z=\calR_i$
% outside of the rectangle.  $\delta$ is the area outside of $A_m$ on
% the tails.

All cases are described in the following context:
\begin{eqnarray*}
  && \calR_i=X^J Y^{1-J} \\
  && g_i \mbox{ drawn from } W=W'^K T^{1-K}
\end{eqnarray*}

\def\Xbar{{\overline X}}

Certain of the cases will also require $V_1$ and $V_2$, defined as
independent random variables distributed
identically to $U_{A_m}$.
The two random variables depend on~$\calR_i$, and hence on~$i$.

\paragraph{Case 0:  ($m\ge\delta$ or $\max_{g\in G}\Pr(\calR_i=g)>1-\delta$)}
Note that $m\ge\delta$ implies $\max_{g\in G}\Pr(\calR_i=g)\ge m\ge\delta$.
Hence,  $\max_{g\in G}\Pr(\calR_i=g)\ge\min(\delta,1-\delta)$ and this
case represents the initial situation,
when the probability distribution of~$G$ still includes at least one
group element whose probability of occurrence is high.
Since $\delta$ is a constant, we need only show that
we can make constant progress.  Specifically, after a constant number
of steps, we need to show that $\max_{g\in
G}\Pr(\calR_i=g)<\min(\delta,1-\delta)$.  Lemma~\ref{lem:maxProb}
shows that if this is true for some $i$, then it will be true for all
$j\ge i$.

One can show for
arbitrary constant~$\delta$ that there large enough constants~$i$
and~$\phi$, such that $|G|\ge\phi$ implies
$\max_{g\in G}\Pr(\calR_i=g)<\min(\delta,1-\delta)$.
We omit the details.
\IGNORE{
{\bf\Large FILL THIS IN,
and make precise ``for large enough~$i$'' in full version of paper.}
}

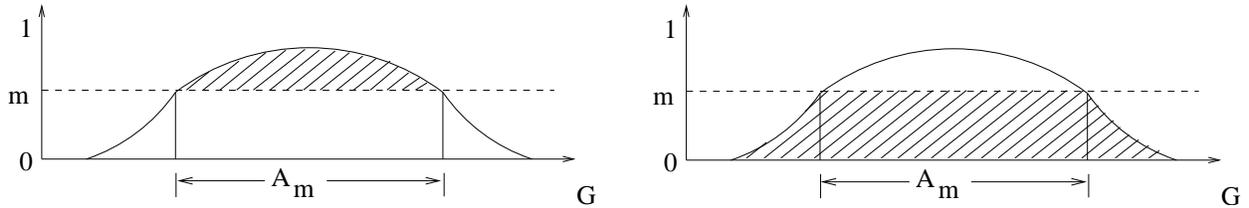
\begin{figure}[htb]
% 70\% reduction from Am-case1[ab].fig
\centerline{\input{Am-case1a.pstex_t}\hspace{0.3truein}
            \input{Am-case1b.pstex_t}}
\caption{Case 1:  Left shaded part is unnormalized probability density
         for~$X$; right shaded part is unnormalized probability
         density for~$W'$ (shaded parts have area less than~1)}
\label{fig:Am-case1}
\end{figure}

\paragraph{Case 1:  ($m<\delta$
         and $\max_{g\in G}\Pr(\calR_i=g)\le1-\delta$ and $m|A_m|<1-\beta$)}

Intuitively, if $\norm{X}^2$ is larger than $\max_{g\in G}\Pr(W'=g)$,
then we will make progress to a more uniform distribution
via Theorem~\ref{thm:reductionBoundedSupport}.  We require that $\norm{X}$ and
$\norm{W}$ be sufficiently large.
We enforce this condition through
$\Pr(\calR_i\notin A_m)\ge\delta$ and through $\Pr(\calR_i\in
A_m)-m|A_m|>(1-\delta)-(1-\beta)=\beta-\delta$.  This allows us to
choose $X$ and $W'$ as in Figure~\ref{fig:Am-case1}.

\begin{sloppypar}
Let $f_1(g)=\max(0, \Pr(\calR_i=g)-m)$.
Let $\Pr(J=1)=\sum_{g\in G}f_1(g)=\sum_{g\in G}f_1(g)=\Pr(\calR_i\in A_m)-m|A_m|$.
Define $X$ so that
$\Pr(X=g)=f_1(g)/\Pr(J=1)=f_1(g)/$\linebreak[0]$\sum_{g'\in G}f_1(g')$.
Let $f_2(g)=\min(m, \Pr(\calR_i=g))$.
Let $\Pr(K=1)=\sum_{g\in
G}f_2(g)=\Pr(\calR_i\notin A_m)+m|A_m|$.
Define $W'$ so that $\Pr(W'=g)=f_2(g)/\Pr(K=1)=f_2(g)/\sum_{g'\in G}f_2(g')$
and $W'$ is independent of~$X$.  Note that
$\Pr(J=1)=\sum_{g\in G}f_1(g)>1-2\delta-(1-\beta)=\beta-2\delta$.
Note that $\Pr(K=1)=\sum_{g\in
G}f_2(g)=\Pr(\calR_i\notin A_m)+m|A_m|$ and hence
$\delta\le\Pr(K=1)<(1-\beta)+2\delta=1+2\delta-\beta$.
\end{sloppypar}

% Note that $\Pr(W'\in{\supp(X)})=1-\delta$.
% So $\Pr(XW'\in{\supp(X)})\le1-\delta$ and
% Theorem~\ref{cor:reductionBySupport}??
% Then use Theorem~\ref{thm:newNextStep}??

We wish to apply Theorem~\ref{thm:reductionBoundedSupport}.
Let $X$ and~$W$ of Theorem~\ref{thm:reductionBoundedSupport}
correspond to $X$ and~$W'$ in our context.  Denote the $\delta$ of
Theorem~\ref{thm:reductionBoundedSupport} by
$\delta'=\delta/\Pr(K=1)>\delta/(1+2\delta-\beta)$ for $\delta$ in our
context.  The conclusion of
the theorem then yields
that for a fixed~$g$ drawn from the distribution of~$W'$, with probability
at least $1-1/\lambda$,
$\norm{Xg^{E_i}}\le\sqrt{\lambda(1-\delta'/2)}\,\norm{X}
<\a\norm{X}$, where
\[
\a = \sqrt{\lambda\frac{2+3\delta-2\beta}{2+4\delta-2\beta}}.
\]

We have $\norm{Xg^{E_i}}/\norm{X}$ bounded above, and we wish to
invoke Theorem~\ref{thm:newNextStep} by identifying $Z$ with~$\calR_i$
and $A=\supp(X)$ with~$A_m$.  The conditions $\Pr(J=0)\Pr(Y=g)=m$ for
$g\in\supp(X)$ and $\Pr(J=0)\Pr(Y=g)\le m$ hold also in our context.
We invoke the theorem with
$\a$ as above,
and $\Pr(Z\notin A_m)=\Pr(\calR_i\notin A_m)
<2\delta$.
Recall that $\Pr(J=1)>\beta-2\delta$.
So,
\begin{eqnarray*}
c&=&(1-\a^2)\*\,
   (\Pr(J=1))^2/(m|A_m|\Pr(Z\notin A_m)+1) \\
%  &=&(1-{\lambda\,(1+(\delta-\beta)/2)/(1+\delta-\beta)})\*\,
%    (\Pr(\calR_i\in A_m)-m|A_m|)^2
%    / (m|A_m|\Pr(\calR_i\notin A_m)+1) \\
%  &\ge&(1-{\lambda\,(1+(\delta-\beta)/2)/(1+\delta-\beta)})\*\,
%    (\Pr(\calR_i\in A_m)-1+\beta)^2
%    / ((1-\beta)\Pr(\calR_i\notin A_m)+1) \\
%  &\ge& (1-{\lambda\,(1+(\delta-\beta)/2)/(1+\delta-\beta)})\*\,
%     (\beta-2\delta)^2 / (1+(1-\beta)(\delta-1/|G|)) \\
 &>& \left(1-\lambda\frac{2+3\delta-2\beta}{2+4\delta-2\beta}\right)\*\,
    \frac{(\beta-2\delta)^2}{1+2(1-\beta)\delta}
  \defeq \cbar \\
 c&<&1-\a^2
\end{eqnarray*}
in Theorem~\ref{thm:newNextStep}.  The random variable~$W$ of
Theorem~\ref{thm:newNextStep} corresponds to $g^{E_i}$ in our current
context
and $Z$ corresponds to $\calR_i$.
To employ Theorem~~\ref{thm:newNextStep}, we also require that $\a<1$,
from which, $c<1-\a^2$ implies $\sqrt{1-c}<1$.  For $\lambda>1$
sufficiently small, $\a<1$.

Hence, $0<\sqrt{1-c}<1$ and $c$ is a constant
determined by $\lambda$, $\delta$ and $\beta$.
So, we have $\norm{\calR_ig^{E_i}}\le\sqrt{1-\cbar}\,\norm{\calR_i}$
with probability at least $(1-1/\lambda)$ for $g$ drawn from the
distribution of~$W'$.  Since $W=W'^KT^{1-K}$, one sees that $\norm{\calR_i
g^{E_i}}\le\sqrt{1-\cbar}\,\norm{\calR_i}$ for $g$ drawn from the
distribution of~$W$ with probability at least
$(1-1/\lambda)\Pr(K=1)\ge(1-1/\lambda)\delta$.

\begin{figure}[htb]
% 70\% reduction from Am-case2[ab].fig
\centerline{\input{Am-case2a.pstex_t}\hspace{0.3truein}
            \input{Am-case2b.pstex_t}}
\caption{Case 2:  Left shaded part is unnormalized probability density
         for~$X$; right shaded part is unnormalized probability
         density for~$W'$
          (shaded part has area less than~1)}
\label{fig:Am-case2}
\end{figure}
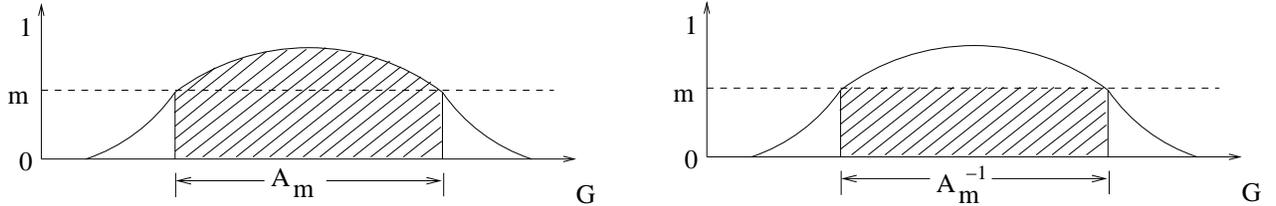

\paragraph{Case 2:  ($m<\delta$
        and $\max_{g\in G}\Pr(\calR_i=g)\le1-\delta$
        and $m|A_m|\ge1-\beta$
        and $\Pr(V_1\inv V_2\in A_m)\le 0.997$)}

Intuitively, if $\Pr(W'X\in A_m)$ is small, then we will make progress
toward a more uniform distribution
via Corollary~\ref{cor:altReductionBySupport}.  We enforce this
through $\Pr(V_1\inv V_2\in A_m)\le 0.997$.  We choose an $X$ close
to~$V_2$ and choose
$W'=V_1\inv$ as in Figure~\ref{fig:Am-case2}.  One knows that
$\norm{X}$ and $\norm{W'}$ are sufficiently large, since
$m|A_m|\ge 1-\beta$.

Let
$X$ be a random variable such that
$\Pr(X=g)=\Pr(\calR_i=g)/\Pr(\calR_i\in A_m)$ for $g\in A_m$ and
$\Pr(X=g)=0$ for $g\notin A_m$.
Set $W'=V_1\inv$.  Set $\Pr(J=1)=\Pr(\calR_i\in A_m)$ and note that
$\Pr(J=1)\ge 1-2\delta$.  Similarly, set $\Pr(K=1)=m|A_m|$ and note that
$\Pr(K=1)\ge 1-\beta$.

\begin{sloppypar}
One wishes to apply Corollary~\ref{cor:altReductionBySupport}
with $X$ and~$W'$.
One shows that $\norm{W'X}\le\sqrt{0.997}\norm{X}$.
By Lemma~\ref{lem:maxProb}, $\Pr(W'X=g)\le\max_{h\in G}\Pr(W'=h)=1/|A_m|$.
So $\norm{W'X}$ is maximized when $\Pr(W'X=g)$ equals $1/|A_m|$ or
equals~$0$ for all~$g$.
\IGNORE{
{\bf\Large WHAT ABOUT FRACTIONAL PROBABILITY FOR
SOME $g'$?}
}
Note that $\supp(X)=A_m$.  Define $Z$ and $\phi=\Pr(W'X\in A_m)$
as in Corollary~\ref{cor:altReductionBySupport}.
Hence, $Z=W'X|(W'X\in A_m)$.  ($Z$ is the random variable~$W'X$
conditioned
on the event $W'X\in A_m$.)
Note that one can write $X\probeq V_2^{J'}Y'^{J'}$ for $\Pr(J'=1)
=m|A_m|/\Pr(\calR_i\in A_m)\ge(1-\beta)/(1-\delta)$.
So $1-\phi=\Pr(W'X\notin A_m)
\ge ((1-\beta)/(1-\delta))\Pr(W'V_2\notin A_m)
\ge 0.003(1-\beta)/(1-\delta)$.
So $\norm{W'X}\le\norm{Z}\le
(1/|A_m|)\sqrt{|A_m|(1-0.003(1-\beta)/(1-\delta))}
\le\sqrt{1-0.003(1-\beta)/(1-\delta)}\,\norm{U_{A_m}}
\le\sqrt{1-0.003(1-\beta)/(1-\delta)}\,\norm{X}
\le\sqrt{0.997}\,\norm{X}$.
\end{sloppypar}

Apply Corollary~\ref{cor:altReductionBySupport}, with $X$ and~$W'$ as above,
and with $c=\sqrt{0.997}$.  With probability
$\Pr(K=1)=m|A_m|\ge1-\beta$, a random~$g$ drawn from the distribution
of~$\calR_i$ is as if $g\inv$ were drawn from the distribution
of~$W'$.
Note that $\phi<1$.
Applying the corollary now yields
$\norm{g^{E_i}X}<\sqrt{\lambda(1+c\phi)/2}\,\norm{X}
< \sqrt{\lambda\,(1+\sqrt{0.997})/2}\,\norm{X}$ with probability at least
$1-1/\lambda$ for $g$ drawn from the distribution of~$W'$.
We require $\lambda>1$ to satisfy $\lambda(1+\sqrt{0.997})/2<1$.

\begin{sloppypar}
We wish to apply Theorem~\ref{thm:nextStep}.
(In fact, a variation of Theorem~\ref{thm:nextStep}
is invoked for $WZ$ instead of for $ZW$.)
The random variable~$Y$ is defined by $\calR_i\probeq X^JY^{1-J}$.
To apply the theorem, we need a positive constant~$c$ such that
$\Pr(J=1)\norm{X}\ge c\Pr(J=0)\norm{Y}$.
Note that
$\Pr(J=1)\norm{X}=\Pr(J=1)\norm{U_{A_m}}$.
Note that
$\Pr(J=0)\norm{Y}\le\sqrt{(2\delta/m)m^2}
 = \sqrt{2m\delta}$.
Recall that $\norm{U_{A_m}}=1/\sqrt{|A_m|}$ by
Lemma~\ref{lem:minNorm-ell2}.
Hence, one can choose
\[
c=(1-2\delta)/\sqrt{2\delta},
\]
since $\Pr(J=1)\norm{X}
\ge \Pr(J=1)\norm{U_{A_m}}
\ge (1-2\delta)\norm{U_{A_m}}
= c\sqrt{2\delta}\, \norm{U_{A_m}}
> c\sqrt{2\delta m|A_m|}\, \norm{U_{A_m}}
= c\sqrt{2m\delta/\norm{U_{A_m}}^2}\, \norm{U_{A_m}}
= c\sqrt{2m\delta}
> c\Pr(J=0)\norm{Y}$.

Theorem~\ref{thm:nextStep} is then invoked with the above~$c$ and with
$\a=\sqrt{\lambda\,(1+\sqrt{0.997})/2}$.  The $W$ and $Z$ of
Theorem~\ref{thm:nextStep} correspond to $g^{E_i}$ and $\calR_i$ in
our context.  So, $\norm{g^{E_i}\calR_i}\le\left(({1+\a
c})/\sqrt{1+c^2}\right) \norm{\calR_i}$ with probability at least
$(1-1/\lambda)$ for $g$ drawn from the distribution of~$W'$.  Since
$W=W'^KT^{1-K}$, one sees that
$\norm{g^{E_i}\calR_i}\le\left(({1+\a c})/\sqrt{1+c^2}\right)
\norm{\calR_i}$ for $g$ drawn from the distribution of~$W$ with
probability at least $(1-1/\lambda)\Pr(K=1)\ge
(1-1/\lambda)\*\,(1-\beta)$.
\end{sloppypar}

\begin{sloppypar}
For the inequality $\norm{g^{E_i}\calR_i}\le\left(({1+\a
c})/\sqrt{1+c^2}\right)\norm{\calR_i}$ to be useful,
we require that $\left(({1+\a c})/\sqrt{1+c^2}\right)<1$.
This is true if $\a<1$ and $c$ is sufficiently large.
For the former, we need only require that $\lambda>1$ be sufficiently
small so that $\a=\sqrt{\lambda\,(1+\sqrt{0.997})/2}<1$.
For the latter, it suffices to make
$\delta$ sufficiently small.
We omit the computation of the explicit requirements for $\delta$.
\end{sloppypar}

\IGNORE{
  As seen from Remark~\ref{rem:nextStepConstants},
$\norm{g^{E_i}\calR_i}\le ((1+\a)/2)\norm{\calR_i}$ for sufficiently
small $\lambda$.  {\bf WHY THIS LAST STATEMENT?}
} % end IGNORE

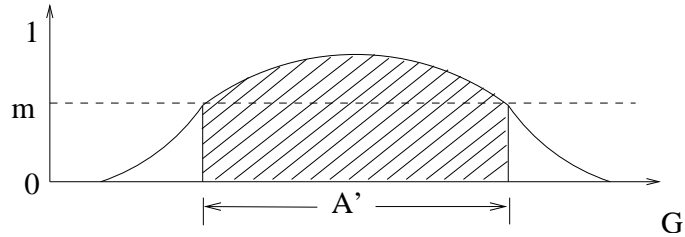
\begin{figure}[htb]
% 80\% reduction from Am-case2[ab].fig
\centerline{\input{Am-case3.pstex_t}}
\caption{Case 3:  Shaded part is unnormalized probability density
         for~$X$
          (shaded part has area less than~1)}
\label{fig:Am-case3}
\end{figure}

\def\Vbar{{\overline V}}

\paragraph{Case 3:  ($m<\delta$
        and $\max_{g\in G}\Pr(\calR_i=g)\le1-\delta$
        and $m|A_m|\ge1-\beta$
        and $\Pr(V_1\inv V_2\in A_m) > 0.997$)}

Intuitively, one constructs an $A'$ close to~$A_m$ with $A'A'$ a
subgroup of~$G$ (Theorem~\ref{thm:escape}).  The argument then splits,
based on whether $A'A'$ is proper in~$G$.  If $A'A'$ is proper in~$G$,
then we choose an $X$ close to $V_1$ as in Figure~\ref{fig:Am-case3}.
The random variable $W=W'$ will be the distribution of random
subproducts on the generators of~$G$.  Under the conditions of Case~3,
one then shows that a random subproduct~$g_i$ drawn from~$W$ has
probability at least 1/2 of satisfying $A'g_i\cap A'=\emptyset$.
Hence, $Xg_i$ escapes from the ``fuzzy subgroup''~$A_m$ with high
probability (Theorem~\ref{thm:escape}).  So $Xg_i^{A_i}$ makes
progress toward a uniform distribution.  If, on the other hand,
$A'A'=G$, then one can show that $U_{A'}\inv U_{A'}$ is already close
to uniform.

We will first construct $A'\subseteq A_m$ such that $A'A'$ is a group.
The random variable~$X$ is then defined such that
$\Pr(X=g)=\Pr(\calR_i=g)/\Pr(\calR_i\in A')$ for $g\in A'$ and
$\Pr(X=g)=0$ for $g\notin A'$.  Let $W=W'$ be the distribution of
random subproducts on the generators of~$G$.  Let
$\Pr(J=1)=\Pr(\calR_i\in A')$.  Note that $\Pr(J=1)\ge1-\delta$.

Since $\Pr(V_1\inv V_2\in A_m)>0.997$, $\Pr(V_1\inv V_2\notin A_m)
= \Pr(V_2\inv V_1=(V_2\inv V_1)\inv\notin A_m\inv)
\le 0.003$.
Recall that $V_1$, $V_2$ and $U_{A_m}$ are identically distributed.
Let $\Vbar_1$ and $\Vbar_2$ be independent random variables with the
same distribution as $U_{A_m\cap A_m\inv}$.
Hence,
$\Pr(\Vbar_1\inv\Vbar_2 \notin A_m\cap A_m\inv)
< \Pr(\Vbar_1\inv\Vbar_2 \notin A_m) + \Pr(\Vbar_1\inv\Vbar_2 \notin A_m\inv)
\le 0.006$.

We claim there exists an $A'\subseteq A_m\cap A_m\inv$ with $A'A'$ a
group,
$|A'A'\setminus A'|\le|A'|/2$,
and $|A'|\ge(9/10)|A_m\cap A_m\inv|$.  To see this, apply
Theorem~\ref{thm:escape} with the constants of
Remark~\ref{rem:constants}.  In particular, $k=20$.  For a random
$(u,v)$ drawn from $\Vbar_1\times \Vbar_2$, $uv\notin A_m\cap A_m\inv$
with probability less than 0.006.  So we conclude from
Theorem~\ref{thm:escape} that there is a $A'\subseteq A_m\cap A_m\inv$
with $A'A'$ a group, and $|(A_m\cap A_m\inv)\setminus A'|<2|A_m\cap A_m\inv|/k
=|A_m\cap A_m\inv|/10$.  So, $|A'|\ge(9/10)|A_m\cap A_m\inv|$.
The inequality $|A'A'\setminus A'|\le|A'|/2$ follows from applying the
constant~$\delta=1/4$ of Remark~\ref{rem:constants} to $|A'A'\setminus
A'|\le\frac{\delta}{1-2\delta}|A'|$ in Theorem~\ref{thm:escape}.

\begin{sloppypar}
We claim $|A_m\cap A_m\inv|> 0.976|A_m|$.  Since $V_1$ and $V_2$ are
independent, $\E(|(A_m\inv V_2)\cap A_m|/|A_m|) = |\{(u,v)\in
A_m\colon~u\inv v\in A_m\}|/|A_m|^2 = \Pr(V_1\inv V_2\in
A_m)>0.997$. Similarly, $\E(|(A_m\inv V_1)\cap A_m\inv|/|A_m|)=\E(|(V_1\inv
A_m)\cap A_m|/|A_m|) =\Pr(V_1\inv V_2\in A_m)>0.997$.  Applying
Lemma~\ref{cor:markovAlt} with its parameter $\lambda=4$ yields
$\Pr(|(A_m\inv V_2)\cap A_m|/|A_m| > 0.988) \ge 3/4$ and $\Pr(|(A_m\inv
V_1)\cap A_m\inv|/|A_m| > 0.988) \ge 3/4.$ So, at least half of the
elements $h\in A_m$ satisfy both $|(A_m\inv h)\cap A_m|/|A_m| > 0.988$ and
$|(A_m\inv h)\cap A_m\inv|/|A_m| > 0.988$.  Choosing one such~$h$ yields
$|A_m\inv\cap A_m|/|A_m| \ge |A_m\inv h\cap A_m\inv\cap A_m|/|A_m| > 0.976$.
\end{sloppypar}

Combining $|A'|\ge(9/10)|A_m\cap A_m\inv|$ and $|A_m\cap A_m\inv|> 0.976|A_m|$
yields $|A'|>0.85|A_m|=0.85|A_m|$.  Recall that $A'\subseteq A_m$, $A'A'$
is a group, and $|A'A'\setminus A'|\le|A'|/2$.  If $|A_m|<(2/3)|G|$,
then $|A'|<(2/3)|G|$ and so $A'A'$ is proper in~$G$ ($A'A'\subset G$).

Assume for the remainder of this case
that $|A_m|<(2/3)|G|$, and hence
$A'A'\subset G$.
We show that
$\norm{Xg^{E_i}}=\norm{X}/\sqrt{2}$ for $g$ drawn from~$W$, with
probability at least~1/2.
Let $W$ be a uniform random variable on the random
subproducts on the generators of the group~$G$.  
By Lemma~\ref{lem:randomSubproduct}, $|A'A'|<|G|$ implies
$\Pr(W\notin A'A')\ge1/2$.
Composing $\Pr(W\notin A'A')\ge1/2$ with
Lemma~\ref{lem:cubeDoubling} implies that
$\Pr(A'W\cap A'=\emptyset)=1/2$.
Since $X=U_{A'}$, $\norm{Xg^{E_i}}=\norm{X}/\sqrt{2}$ with probability
at least~1/2 for $g$ drawn from~$W$.

\begin{sloppypar}
We wish to apply Theorem~\ref{thm:nextStep}.
The random variable~$Y$ is defined by $\calR_i\probeq X^JY^{1-J}$.
To apply the theorem, we need a positive constant~$c$ such that
$\Pr(J=1)\norm{X}\ge c\Pr(J=0)\norm{Y}$.
Recall that $\norm{U_{A_m}}=1/\sqrt{|A_m|}$ by
Lemma~\ref{lem:minNorm-ell2}, and similarly $\norm{U_{A'}}=1/\sqrt{|A'|}$.
Note that $|A'|>0.85|A_m|$ and $\Pr(J=1)\ge 1-\beta$ implies
$\Pr(J=1)\norm{X}=\Pr(J=1)\norm{U_{A'}}
= \sqrt{|A'|/|A_m|}\*\,\Pr(J=1)\norm{U_{A_m}}
> 0.85(1-\beta)\norm{U_{A_m}}$.
Note that
$\Pr(J=0)\norm{Y}\le\sqrt{(2\delta/m)m^2}
= \sqrt{2m\delta}
< \sqrt{2m\delta}$.
Hence, one can choose
\[
c = 0.85(1-\beta)/\sqrt{2\delta},
\]
since $\Pr(J=1)\norm{X}
> 0.85(1-\beta)\norm{U_{A_m}}
= c\sqrt{2\delta}\, \norm{U_{A_m}}
> c\sqrt{2\delta m|A_m|}\, \norm{U_{A_m}}
= c\sqrt{2m\delta/\norm{U_{A_m}}^2}\, \norm{U_{A_m}}
= c\sqrt{2m\delta}
\ge c\Pr(J=0)\norm{Y}$.
\end{sloppypar}

Theorem~\ref{thm:nextStep} is then invoked with the above~$c$ and with
$\a=1/\sqrt{2}$.  The $W$ and $Z$ of
Theorem~\ref{thm:nextStep} correspond to $g^{E_i}$ and $\calR_i$ in
our context.  So, $\norm{g^{E_i}\calR_i}\le\left(({1+\a
c})/\sqrt{1+c^2}\right) \norm{\calR_i}$ with probability at least
$1-1/\lambda$ for $g$ drawn from the distribution of~$W$.

\begin{sloppypar}
For the inequality $\norm{g^{E_i}\calR_i}\le\left(({1+\a
c})/\sqrt{1+c^2}\right)\norm{\calR_i}$ to be useful,
we require that $\left(({1+\a c})/\sqrt{1+c^2}\right)<1$.
This is true if $c>\sqrt{2}$.  For this, it suffices to make
$1-\beta>2\sqrt{\delta}$.
\end{sloppypar}

\bigskip
The preceding analysis demonstrates the following lemma.

\def\Jbar{{\bar J}}

\begin{lemma}
\label{lem:main}
Let $\overline\calR_i$ be independent and identically distributed to~$\calR_i$.
For any choice of positive
parameters $a$, $b$ and~$c$ in the Fibonacci cube algorithm, there are
constants $\alpha<1$, $\rho>0$, $\beta>0$ and $\iota>0$
such that for $i>\iota\log|G|$ one of the following holds:
\begin{enumerate}
\item[(i)] $\norm{\calR_{i+1}}\le \alpha\norm{\calR_i}$ with probability at
least~$\rho$; or
\item[(ii)] $\Pr(\calR_i\inv\overline\calR_i=g)\ge(3/4)(1-\beta)^2/|G|$
for all $g\in G$.
\end{enumerate}
Further, Let $\phi>1$.
For $i\ge\iota\log|G|+(\phi/\rho)\,(1+(1/2)\log_{1/\alpha}|G|)$,
case~ii above occurs with probability at least
$1-\exp(-(\phi(1-1/\phi)^2/4)\,\log_{1/\alpha}|G|)$.
\end{lemma}
\begin{proof}
The proof follows from the analysis of the three cases just presented.
As discussed in the analysis of Case~0, after a constant number of
steps of the Fibonacci cube algorithm, Case~0 will never again be
revisited, with high probability.  Therefore, after $\iota\log|G|$ steps,
for some constant~$\iota$, the probability of ever revisiting Case~0 will
be less than $\exp(-\log|G|)$.  Hence, we can ignore Case~0 for purposes
of the analysis.

We show that $|A_m|<(2/3)|G|$ implies case~i and that
$|A_m|\ge(2/3)|G|$ implies case~ii.
Assume first that $|A_m|<(2/3)|G|$.
In each of the three cases, we concluded that
$\norm{\calR_{i+1}}=\norm{\calR_i g^{E_i}}\le \alpha\norm{\calR_i}$
or $\norm{\calR_{i+1}}=\norm{g^{E_i}\calR_i}\le \alpha\norm{\calR_i}$
with probability at least $\rho$ for appropriate $\alpha<1$ and $\rho>0$.
(In Case~3, this conclusion need not hold if $|A_m|\ge(2/3)|G|$.)
The parameters~$\alpha$ and~$\rho$ are defined in terms of
$\beta$, $\delta$, $\lambda$ and~$G$ for each of the three cases.

In order to make the parameters~$\alpha$
independent of the particular case, one chooses $\alpha$ to
be the maximum of the three definitions for each of the three cases.
In order to make $\rho$ independent of the particular
case, define $\rho_1$, $\rho_2$ and~$\rho_3$ to be the probabilities
for the three cases.  Then let $\rho=\min(\rho_1/(ad),\rho_2/(bd),\rho_3/(cd))$
for $d=1/a+1/b+1/c$.  In particular, $\rho$ can be maximized by choosing
$a=\rho_1$, $b=\rho_2$ and $c=\rho_3$, whereupon
$\rho=1/(1/\rho_1+1/\rho_2+1/\rho_3)$.

It remains to verify that the constants $\beta$, $\delta$ and~$\lambda$ can be
simultaneously chosen to meet the requirements of the analysis in
Cases~1, 2 and~3.
Recall that $1>\beta>2\delta>0$ and $\lambda>1$.
Collecting the bounds from Case~1, we require $\lambda>1$ to be
sufficiently small that
$\alpha=\sqrt{\lambda\frac{2+3\delta-2\beta}{2+4\delta-2\beta}}$.
Collecting the bounds from Case~2, we require that $\lambda>1$ such
that $\lambda(1+0.997)/2<1$.  We further require that $\delta$ be
sufficiently small
to satisfy $(1+\alpha c)/\sqrt{1+c^2}<1$ for
$c=({1-2\delta})/\sqrt{2\delta}$.
The bounds from Case~3 require that $\lambda>1$ and
$1-\beta>2\sqrt{\delta}$.

There can be at most $\log_{1/\alpha}\sqrt{|G|}=O(\log|G|)$ distinct instances
of~$i$ such that $\norm{\calR_{i+1}}>\alpha\norm{\calR_i}$.  To see this,
note that $\norm{\calR_0}=1$ and
$\norm{\calR_i}\ge\norm{U_G}=1/\sqrt{|G|}$ for all~$i$ by
Lemma~\ref{lem:minNorm-ell2} and that $\norm{\calR_{i+1}}\le\norm{\calR_i}$
by Lemma~\ref{lem:bounded-ell2}.

% TECHNICALLY, NEED TO WORRY ABOUT PRODUCING RIGHT $A_m$ at one step,
% and then wrong $A_m$ again at the next step.  Must choose right~$t$.

With the probability in the statement of the lemma, we must show we
are in Case~3 and $|A_m|\ge(2/3)|G|$ with the
stated probability after the stated number of steps.
We will then show that this implies case~ii.
We define the $i$-th step to be a {\em success} if
$\norm{\calR_{i+1}}\le\alpha\norm{\calR_i}$.
So, at most $\log_{1/\alpha}\sqrt{|G|}$
successes may occur
for distinct~$i$.
We know that for a given~$i$, a success will occur with probability at
least $\rho$, or else
$|A_m|\ge(2/3)|G|$.

Consider Chernoff's bound (Theorem~\ref{thm:chernoff}).
Assume a success with probability at most~$p=\rho$,
and assume $t=(1+\log_{1/\alpha}\sqrt{|G|})/(\rho(1-\epsilon))$ trials.
Chernoff's bound predicts at least $\lfloor (1-\epsilon)pt
\rfloor
\ge\log_{1/\alpha}\sqrt{|G|}$ successes over $t$~trials
with probability at least
$1-\exp(-{\epsilon}^2 pt/2)$.
We have seen that more than $\log_{1/\alpha}\sqrt{|G|}$ successes are
impossible.
So, with probability at least $1-\exp(-{\epsilon}^2 pt/2)$,
we are in Case~3 and $|A_m|\ge(2/3)|G|$ for some step~$j$ among the
first $t$ steps.  Let $\epsilon=1-1/\phi$ for $\phi>1$.  This yields
the probability of the lemma.

Hence, there is a~$j$ such that $\calR_j$ is in Case~3 and $|A_m|\ge(2/3)|G|$.
Combining the condition $m|A_m|\ge1-\beta$ of Case~3 with
$|A_m|\ge(2/3)|G|$ implies that $m\ge({1-\beta})/((2/3)|G|)$.  Define a
$\{0,1\}$-random variable~$J$ such that $\Pr(J=1)=m|A_m|$.  Note
$\Pr(J=1)\ge1-\beta$.  Let $\calR_j\probeq X^J Y^{1-J}$ for $X=U_{A_m}$.
Let $\Jbar$ be independent and
distributed identically to~$J$.  Similarly, let $\overline X$ be
independent and distributed identically to~$X$.
Then for $V$ an arbitrary $G$-valued random variable,
$\calR_j\inv V\overline\calR_j\probeq(X\inv V\overline X)^{J\Jbar}
Y'^{1-J\Jbar}$ for some $G$-valued independent random
variable~$Y'$.

\begin{sloppypar}
We show that for an arbitrary $G$-valued random variable~$V$,
$\Pr(\calR_j\inv V\overline\calR_j=g)\ge(3/4)({1-\beta})^2/|G|$
for all $g\in G$ when $|A_m|\ge(2/3)|G|$.  With probability at least
$({1-\beta})^2$,
$J=\Jbar=1$.  Hence, with probability at least $({1-\beta})^2$,
we can take $\calR_j\inv V\overline\calR_j=X\inv V\overline X$.
Applying Lemma~\ref{lem:cubeUniform} with $A=A_m$ and $\alpha=2/3$,
one sees that $X\inv V\overline X$ is 1/2-uniform and that
$\Pr(X\inv V\overline X=g)\ge(|G|/3)/|A_m|^2\ge(3/4)/|G|$.
\end{sloppypar}

% Next, $\Pr(X\inv\overline X=g)\ge (3/4)/|G|$.
% To see this, note that
% \begin{eqnarray*}
% \lefteqn{|\{(h_1,h_2)\colon~h_1\inv h_2=g\mbox{ and
%               }h_1,h_2\in A_m\}|}\qquad \\
% &=& |\{(h_1,h_2)\colon~h_1\inv h_2=g\}| 
%  - |\{(h_1,h_2)\colon~h_1\inv h_2=g\mbox{ and }
%   h_1\notin A_m\mbox{ or }h_2\notin A_m\}| \\
% &\ge& |G|-2|G\setminus A_m| \\
% &\ge& |G|/3.
% \end{eqnarray*}

% We show $\Pr(\calR_j\inv\overline\calR_j=g)\ge(3/4)({1-\beta})^2/|G|$
% for all $g\in G$ when $|A_m|\ge(2/3)|G|$.  With probability at least
% $({1-\beta})^2$,
% $J=\Jbar=1$.  Hence, with probability at least $({1-\beta})^2$,
% we can take $\calR_j\inv\overline\calR_j=X\inv\overline X$.
% Next, $\Pr(X\inv\overline X=g)\ge (3/4)/|G|$.
% To see this, note that
% \begin{eqnarray*}
% \lefteqn{|\{(h_1,h_2)\colon~h_1\inv h_2=g\mbox{ and
%               }h_1,h_2\in A_m\}|}\qquad \\
% &=& |\{(h_1,h_2)\colon~h_1\inv h_2=g\}| 
%  - |\{(h_1,h_2)\colon~h_1\inv h_2=g\mbox{ and }
%   h_1\notin A_m\mbox{ or }h_2\notin A_m\}| \\
% &\ge& |G|-2|G\setminus A_m| \\
% &\ge& |G|/3.
% \end{eqnarray*}
% So, $\Pr(X\inv\overline X=g)\ge(|G|/3)/|A_m|^2\ge(3/4)/|G|$.

We have seen
$\Pr(\calR_j\inv V\overline\calR_j=g)\ge(3/4)({1-\beta})^2/|G|$.
We show that
$\Pr(\calR_i\inv\overline\calR_i=g)\ge(3/4)\*({1-\beta})^2/|G|$
for all $i\ge j$.  To see this, define
$\calX=\calR_i\inv\overline\calR_i
=V_1\calR_j\inv V_2\overline\calR_j V_3$.
For $U$ uniform on~$G$, we can write
$\calR_j\inv V_2\overline\calR_j\probeq U^JY^{1-J}$ for $J$ an independent
$\{0,1\}$-random
variable with $\Pr(J=1)=(3/4)({1-\beta})^2$.
So $\calX\probeq
\left(V_1 U V_3\right)^J
 \left(V_1 Y V_3\right)^{1-J}$.
Lemma~\ref{lem:uniform} shows that $V_1 U V_3$ is uniform.
So $\Pr(\calX=g)\ge\Pr(J=1)/|G|=(3/4)({1-\beta})^2/|G|$ for all $g\in
G$.
\end{proof}

For some applications, Lemma~\ref{lem:main} may suffice, since it
promises to produce each group element with a minimum probability
$(3/4)(1-\beta)^2/|G|$.  For an $\varepsilon$-uniform random distribution,
one must do a little more.  The next section is concerned with producing an
$\varepsilon$-uniform distribution.

\section{Constructing $\varepsilon$-uniform from $\varepsilon$-semi-uniform}
\label{sec:semi-to-uniform}

Lemma~\ref{lem:main} shows that for $i$ sufficiently large, Algorithm
Fibonacci Cube constructs an $\alpha$-semi-uniform random variable,
$\calR_i\inv\overline\calR_i$, with the stated probability for
$\alpha=(3/4)({1-\beta})^2$.  This section shows that constructing a
$\varepsilon$-semi-uniform random distribution is tantamount to
constructing a $\varepsilon$-uniform random distribution.  This is
shown in the next theorem uses $W=\calR_i\inv\overline\calR_i$ in
order to efficiently construct a $\beta$-uniform random variable.

\begin{theorem}
\label{thm:semiuniform}
Let $G$ be a group.  Let $W$ be an $\alpha$-semi-uniform random
variable on~$G$.  Let $\calP_0$ be an arbitrary $G$-valued random variable.
Let $E_i$ be independent, uniform random variables on~$\{0,1\}$.
For all $i>0$, define
$\calP_{i+1}=\calP_i g_i^{E_i}$ for $g_i$ drawn from the distribution
of $W$.
Let $\gamma=14/(11+3\alpha)$.
Then, $\Pr(\calP_t=g)\le7/8$ for
$t\ge2\log_\gamma|G|+\log_\gamma(64\lambda)$,
with probability at least $1-1/\lambda$.
Hence $W\calP_t$ is a $\max(\alpha,7/8)$-uniform random variable
with probability at least $1-1/\lambda$.  
\end{theorem}

\begin{proof}
Define the set $A_i=\{g\colon~
\Pr(\calP_i=g)\ge (7/4)/|G|\}$.
Note that $|G\setminus A_i|\ge (3/7)|G|$, since
otherwise $|A_i|>(4/7)|G|$, which implies $\Pr(\calP_i\in
A_i)=\sum_{g\in A_i}\Pr(\calP_i=g)
>|A_i|(7/4)/|G|>1$.

Define $T_i=\sum_{h\in A_i}(\Pr(\calP_i=h)-(7/4)/|G|)^2$ for $i\ge0$.
We will find an upper bound on $\E(T_{i+1})$ as compared to~$T_i$.
Define
$x_h=\Pr(\calP_i=h)-(7/4)/|G|$.  Hence
$x_{h(g\inv)^{E_i}}=x_h/2+x_{hg\inv}/2$ since $E_i$ and $\calP_i$ are
independent.  Note that for $i\ge 0$,
\[
T_i=\sum_{h\in A_i}(\Pr(\calP_i=h)-(7/4)/|G|)^2 
=\sum_{h\in G}\left(\max\left(0,x_h\right)\right)^2.
\]

We show that $T_{i+1}\le T_i$ for any value of~$g_i$.
Recall that $\calP_{i+1}=\calP_i g_i^{E_i}$.
\begin{eqnarray*}
T_{i+1}&=& \sum_{h\in A_{i+1}}
        \left(\Pr(\calP_ig_i^{E_i}=h)-(7/4)/|G|\right)^2 \\
  &=& \sum_{h\in G}\left(\max(0, x_{h(g_i\inv)^{E_i}})\right)^2 \\
  &=& \sum_{h\in G}\left(\max(0, x_h/2 + x_{hg_i\inv}/2)\right)^2 \\
  &\le& \frac{1}{4} \sum_{h\in G}\Bigl(\max(0, x_h)\Bigr)^2
        + \frac{1}{4}\sum_{h\in G}\left(\max(0, x_{hg_i\inv})\right)^2
        + \frac{1}{2} \sum_{h\in G}\max(0, x_h)\max(0, x_{hg_i\inv}) \\
  &\le& T_i/4 + T_i/4 + T_i/2 \\
  &=& T_i
\end{eqnarray*}
where the Cauchy-Schwartz inequality was invoked to show
\[
\sum_{h\in G}\max(0, x_h)\max(0, x_{hg_i\inv})
\le\sqrt{\sum_{h\in G}(\max(0, x_h))^2}\,\sqrt{\sum_{h\in G}(\max(0,
x_{hg_i\inv}))^2} = T_i.
\]

% Define the set $B_i=\{g\colon~
% \Pr(\calP_i=g)\le (5/4)/|G|\}$.
% Note that $|B_i|\ge |G|/5$, since
% otherwise $|G\setminus B_i|>(4/5)|G|$, which implies
% $\Pr(\calP_i\notin B_i)=\sum_{g\notin B_i}\Pr(\calP_i=g)
% >|G\setminus B_i|(5/4)/|G|>1$.

Since $W$ is $\alpha$-semi-uniform,
by Lemma~\ref{lem:decomposition} we can write
$W=U^J V^{1-J}$ for $G$-valued random variables $U$ and $V$, with
$U$, $V$ and $J$ independent, $U$~uniform, and $Pr(J=1)=1-\alpha$.
Note that $2(\frac{x_h}{2}+\frac{x_g}{2})^2 \le
(x_h^2+x_g^2)$ follows from elementary algebra.  Note that $x_g\le
0$ for $g\notin A_i$.  The notation $E_{g_i\in U}(f(g_i))$ denotes
$\E(f(U))$ for a function~$f(\cdot)$ from $G$ to the real numbers.
Since $U$ and $\calP_i$ are independent, if one conditions on
$J=1$ (implying that $g_i$ is drawn from~$U$), then the
following is true.
\begin{eqnarray*}
\lefteqn{\E_{g_i\in U}(T_{i+1}\mid J=1)}\quad && \\
  &=& \E_{g_i\in U}\left(\sum_{h\in A_{i+1}}
            \left(\Pr(\calP_ig_i^{E_i}=h)-(7/4)/|G|\right)^2\right) \\
  &=& \E_{g_i\in U}\left(\sum_{h\in G}
            \left(\max\left(0,x_{h(g_i\inv)^{E_i}}\right)\right)^2\right) \\
  &=& \frac{1}{|G|} \sum_{h\in G}\sum_{g\in G}
                \left(\max(0,\frac{x_h}{2}+\frac{x_{hg\inv}}{2})\right)^2 \\
  &=& \frac{1}{|G|} \sum_{h\in G}\sum_{g\in G}
                \left(\max(0,\frac{x_h}{2}+\frac{x_g}{2})\right)^2 \\
  &=& \frac{1}{|G|} \sum_{\begin{array}{c}
        \scriptstyle{h\in A_i}\\\scriptstyle{g\in A_i} \end{array}}
                \left(\max(0,\frac{x_h}{2}+\frac{x_g}{2})\right)^2
    + \frac{2}{|G|} \sum_{\begin{array}{c}
        \scriptstyle{h\in A_i}\\\scriptstyle{g\notin A_i}
        \end{array}}
                \left(\max(0,\frac{x_h}{2}+\frac{x_g}{2})\right)^2 \\
  && \mbox{}  +  \frac{1}{|G|} \sum_{\begin{array}{c}
        \scriptstyle{h\notin A_i}\\\scriptstyle{g\notin A_i} \end{array}}
                \left(\max(0,\frac{x_h}{2}+\frac{x_g}{2})\right)^2 \\
  &\le& \frac{|A_i|}{|G|} \sum_{h\in A_i} x_h^2
   + \frac{2|G\setminus A_i|}{|G|} \sum_{h\in A_i} x_h^2/4 \\
  &=& \frac{|A_i|}{|G|} T_i + \frac{|G\setminus A_i|}{2|G|} T_i
\end{eqnarray*}

Recalling that $|G\setminus A_i|\ge (3/7)|G|$, one sees
\[
\E_{g_i\in U}(T_{i+1}\mid J=1)
 \le \frac{|A_i|}{|G|} T_i + \frac{|G\setminus A_i|}{2|G|} T_i
 \le \frac{11}{14} T_i.
\]

Let $\gamma=14/(11+3\alpha)$.
Then $\E(T_{i+1})\le\E(T_i)/\beta$.  To see this, note
$\E(T_{i+1})=\Pr(J=1)\E(T_{i+1}\mid J=1)+\Pr(J=0)\E(T_{i+1}\mid
J=0)
\le(1-\alpha)(11/14)\E(T_i)+\alpha\E(T_i)=(11/14+3\alpha/14)\E(T_i)$.
An easy argument implies $\E(T_{i+k}\le\E(T_i)/\beta^k$.

Let $\lambda>1$ and let
$t\ge2\log_\gamma(8\sqrt{\lambda}|G|)
 =\log_{1/\gamma}(1/(8\sqrt{\lambda}|G|)^2)$.
Since $T_0\le 1$, $\E(T_t)\le 1/(8\sqrt{\lambda}|G|)^2$.
For $\lambda>1$ in Markov's inequality
(Lemma~\ref{lem:markov}),
one has
$\Pr(T_t< 1/(8|G|)^2)
= \Pr(T_t< \lambda/(8\sqrt{\lambda}|G|)^2)
\ge \Pr(T_t<\lambda\E(T_t))
\ge 1-1/\lambda$. So,
\[
  \Pr(T_t< 1/(8|G|)^2) \ge 1-1/\lambda 
 \mbox{ for }t\ge2\log_\gamma(8\sqrt{\lambda}|G|).
\]

Note that $T_t=\sum_{h\in A_t}(\Pr(\calP_t=h)-(7/4)/|G|)^2 \le
1/(8|G|)^2$
implies that $\max_{h\in G}(\Pr(\calP_t=h)-(7/4)/|G|)^2 \le 1/(8|G|)^2$.
So $\max_{h\in G}\Pr(\calP_t=h) \le 15/(8|G|)$.
If $\Pr(\calP_t=g)\le(15/8)/|G|$ for all $g\in G$, then
by Lemma~\ref{lem:maxProb}, $(1-\alpha)/|G|\le\min_{g\in
G}\Pr(W=g)\le\max_{g\in G}\Pr(W\calP_t=g)\le\max_{g\in
G}\Pr(\calP_t=g)<(15/8)|G|$.  Hence
$\calP_t$ is $\max(\alpha,7/8)$-uniform with the given probability.
\end{proof}

\begin{corollary}
\label{cor:semiuniform}
Assume a random variable~$X$
on~$G$ is $\alpha$-semi-uniform.
Assume it costs $c$~group operations to compute a group element drawn
from the distribution of~$X$.
There is a fixed constant~$\gamma$ such that
one can construct a $\gamma$-uniform random variable~$Y$
for which one can draw a group element from the distribution of~$Y$
using
$O(c+\log|G|/({1-\alpha}))$ group operations.
The cost of constructing~$Y$ is $O(c\log|G|/({1-\alpha}))$ group
operations.
\IGNORE{
{\bf What about dependence on~$\gamma$?  Also, split into two steps:
first find constant-uniform, and then use doubling to get more
uniform.}
}
\end{corollary}

The proof of the corollary is clear.

\begin{theorem}
\label{thm:main}
Let $G=\gen{S}$ be a black box group with $|G|<L$.  One can construct
a $\varepsilon$-uniform~$X$ such that the cost of computing a group
element from the distribution of~$X$ is
$O((\log(1/\varepsilon))\log|G|)$ group operations.
The cost of constructing~$X$ is $O(\log^2|G|+|S|\log|G|)$.
Where $|G|$ is not known a~priori, one can replace $|G|$ by~$L$ in the
asymptotic estimates.
\end{theorem}

\begin{proof}
The timing of the Fibonacci Cube algorithm is immediate since
$O(\log|G|+|S|)$ group
operations are required to compute each $g_i$, and
$t=O((\log1/\varepsilon)\log|G|)$.  So, the
timing of the pseudo-code is $O(\log^2|G|+|S|\log|G|)$.  To compute an
element from the distribution of
$\calR_t\inv\overline\calR_t$ then requires
$t=O((\log1/\varepsilon)\log|G|)$ group multiplications, where each
factor $g_i^{E_i}$ of $\calR_t$ contributes at most one to the number of
multiplications.

Lemma~\ref{lem:main} shows that one can construct an
$\alpha$-semiuniform random variable~$X_1=\calR_t\inv\overline\calR_t$
for $\alpha=(3/4)(1-\beta)^2$ in
$O(b\log|G|)$ steps for $b=(\phi/\rho)/\log(1/\alpha)$.
Hence, $O(b^2\log^2|G|+|S|\log|G|)$ group operations are required to
construct~$X_1$.
Computing an element from the distribution of~$X_1$ costs
$O(b\log|G|)$ group operations.

\begin{sloppypar}
Corollary~\ref{cor:semiuniform} shows that one can construct a
$\gamma$-uniform random variable~$X_2$ using
$O(b\log^2|G|/({1-\alpha})+|S|\log|G|)$
group operations.  One can compute a group element from the distribution
of~$X_2$ using $O((b+1/({1-\alpha}))\log|G|)$ group operations.
\end{sloppypar}

Since $\beta$, $\alpha$, $\phi$ and $\rho$ are all constants,
this implies that one can construct $X_2$ using $O(\log^2|G|+|S|\log|G|)$ group
elements and one can compute an element from~$X_2$ using $O(\log|G|)$
group elements.

It remains to construct an $\varepsilon$-uniform random element from
the give  $\gamma$-uniform random element for arbitrary $\varepsilon>0$.
We take the product of $\lceil\log_2\varepsilon/\log_2\gamma\rceil$ many
$\gamma$-uniform elements drawn from the distribution
of~$\calR_t\inv\overline\calR_t$.
By Lemma~\ref{lem:acceleratorUniform}, this
suffices to produce an $\varepsilon$-uniform random element.

To compute an $\varepsilon$-uniform random element
requires $O(\log1/\varepsilon)$ $\gamma$-uniform random elements.  So,
the number of operations to produce an $\varepsilon$-uniform random
element is $O((\log1/\varepsilon)\log|G|)$.
\end{proof}

\IGNORE{

{\bf
\begin{sloppypar}
We also have to say how small $\ell^2$ norm converts to small
$\varepsilon$ in $\varepsilon$-uniform.  It's easy to see that $X$ being
$\varepsilon$-uniform implies that $\norm{X}\le
\sqrt{|G|}(1/|G|+\varepsilon)\le 1/\sqrt{|G|}+\varepsilon/\sqrt{G}$.
We have to argue that if
$\norm{X}=r+1/\sqrt{|G|}$, then $\max_{g\in G}\Pr(X=g)$ is maximized
when $\Pr(X=g')=p_r>1/\sqrt{|G|}$ and $\Pr(X=g)=k_r<1/\sqrt{|G|}$ for all
$g\not= g'$.  Then $X$ will be $p-1/\sqrt{|G|}$-uniform or
$p-k$-uniform since $k<1/\sqrt{|G|}$.  So
$r+1/\sqrt{|G|}=\norm{X}=\sqrt{(p-k)^2+2(p-k)k+|G|k^2}$.

Since $p_r+k_r(|G|-1)=1$, $k=(1-p)/(|G|-1)\ge(1-p)/|G|$.
So,
$r+1/\sqrt{|G|}=\norm{X}\ge\sqrt{(p-k)^2+2(p-k)(1-p)/|G|+(1-p)^2/|G|}
\ge\sqrt{(p-k+(1-p)/|G|)^2+(1-p)^2(1/|G|-1/|G|^2)}
$.  {\bf FINISH THIS!!!}
\end{sloppypar}
}

{\bf This is converting from $\ell^2$ to $\ell^\infty$ norm.  Minkowski
inequality?}

}

\begin{remark}
Chernoff's bound shows that the probability of error can be further
reduced by a power of~$n$ at the cost of multiplying $t$ by the
factor~$n$.
\IGNORE{
Given the parameters $\alpha$ and $\beta$ from this section, one can
choose the parameters $a$, $b$ and~$c$ so as to produce a minimal
coefficient of complexity in the bound $O(\log^2|G|)$.  However, the
bound derives from a worst case analysis over all possible groups~$G$
and may not represent typical groups.  Hence, some applications may
select $a$, $b$ and~$c$ for a minimal coefficient of complexity, while
other applications may choose $a$, $b$ and~$c$ heuristically to
perform well for ``typical'' groups.  The asymptotic complexity bound
holds in both cases.
}% end IGNORE
\end{remark}

Theorem~\ref{thm:main} states a complexity of~$O(\log^2|G|+|S|\log|G|)$ group
operations.  In the unusual case that $|S|>O(\log|G|)$, there is a
black box algorithm to quickly produce a smaller generating
set~\cite{BabaiCoopermanetal91,CoopermanFinkelstein93}.
We quote that theorem here.
\begin{theorem}[{from~\cite[Theorem~2.3]{BabaiCoopermanetal91}}]
Let $G=\gen{S}$ be a finite group.  Let $L$ be a known upper bound on the
length of all subgroup chains in~$G$.  Then for any fixed parameter~$p$ such
that $0<p<1$, with probability at least~$p$ one can find a generating set~$S'$
with $|S'|=O(L\log(1/(1-p)))$, using $O(|S|\log L\log(1/(1-p)))$ group
operations.
\end{theorem}

\IGNORE{
\begin{remark}
It follows from Erd\H os and
R\'enyi~\cite[Theorem~1]{ErdosRenyi65} that once one has accumulated
$k\ge 2\log|G|$\linebreak[0]$\varepsilon$-uniform random group
elements (the exact
number depends on $\varepsilon$ and the desired probability of error),
one can replace the cube $g_1^{E_1}\cdots g_t^{E_t}$ by a cube of
length~$k$, which is usually shorter.  Babai~\cite{Babai91} used this
approach to produce additional group elements in $O(\log|G|)$ time.
Although the Fibonacci cube algorithm already produces additional
elements in $O(\log|G|)$, the constant~$k$ is likely to be smaller
than~$t$, thus saving a constant factor in time.
\end{remark}
} % end IGNORE

% end IGNORE

\section{Experimental Results}
\label{sec:experiments}

The current results are highly preliminary.  For the Fibonacci cube
algorithm, we initialize the first elements of the cube to be the
group generators.  We take the parameters $a=b=c=1$ as a simple
heuristic choice.  We compute only $\calR_t\inv\overline\calR_t$,
which, in principle, is $\varepsilon$-semi-uniform, but not
necessarily $\varepsilon$-uniform.  We take $t=20$, 25, and~30.  After
the precomputation of the $g_1,\ldots,g_{20}$ that determine
$\calR_t\inv\overline\calR_t$, we draw 10,000 elements from the
distribution of $\calR_t\inv\overline\calR_t$.

The table shows the results of tests on the distribution of
$\calR_t\inv\overline\calR_t$ according to a partition into conjugacy
classes.  (The conjugacy class of $g\in G$ is $\{g^h\colon h\in G\}$.)
The groups tested on are all simple groups.  Later experiments will
consider other parameters than $a=b=c=1$.  They will incorporate the
ideas of Section~\ref{sec:semi-to-uniform}.  They will also look at
distributions over other group partitions than that of conjugacy
classes.

The $\chi^2$ distribution was applied with a critical value of~0.05.
The $\chi^2$ test accepts the hypothesis of uniform randomness when
the observed $\chi^2$ statistic satisfies
$\chi^2<\chi^2_{\lower1pt\hbox{$\scriptstyle .05$}}$.

The number of degrees of freedom in the $\chi^2$ test is one less than
the number of conjugacy classes.  However, in most tests, the smaller
conjugacy classes showed fewer than five observations.  Hence, the
smallest conjugacy classes have been merged so that the smallest set
in the partition has just enough conjugacy classes to have at least
five observations.  The number of degrees of freedom is then adjusted
accordingly, as one less than the number of final categories.

These experimental results are intended only to demonstrate the
quality of the random elements in computer experiments.  In principle,
the distribution ${g_{20}\inv}^{\overline
E_{20}},\ldots,{g_1\inv}^{\overline
E_20},g_1^{E_1},\ldots,g_{20}^{E_{20}}$ can produce at most
$2^{40}\approx 10^{12}$ group elements.  This is not too much larger
than the order of the groups being tested.  Hence, the experimental
distribution of the individual group elements is most likely {\em not}
close to uniform.  However, the $\chi^2$ test shows that an empirical
computation will not be able to distinguish the distribution of group
elements according to conjugacy class from a distribution based on
uniformly random group elements.

\begin{table}
\begin{center}
% \begin{tabular}{|c|p{4.5in}|p{1.2in}|}
\begin{tabular}{|c|c|c|c|c|c|}
\hline
Group~$G$ & $|G|$ & terms ($t$)/ & total num classes/
     & $\chi^2$ probability & $\chi^2_{\lower1pt\hbox{$\scriptstyle .05$}}$ \\
        &     & precomp.{} (group op's)     & $\chi^2$ degrees freedom
     &                      &       \\
\hline
$M_{24}$ & $2.4\times 10^8$           & 20/60  & 26/22  & 12.3  & 33.9 \\
$McL$ & $9.0\times 10^8$              & 25/98  & 24/20  & 35.6 & 31.4 \\
${\rm SL}(7,2)$ & $1.6\times 10^{14}$ & 25/110 & 117/98 & 112.0 & 122.1 \\
$Suz$ & $4.5\times 10^{11}$           & 30/184 & 43/32  & 27.3 & 46.2 \\
$A_{15}$ & $6.5\times 10^{11}$        & 30/204 & 94/68  & 52.0 & 88.2 \\
\hline
\end{tabular}
\end{center}
\end{table}

The $\chi^2$ test accepts the hypothesis at the 0.05 significance
level for all groups, except the McLaughlin group (McL).  The
McLaughlin group is accepted at the 0.01 significance level.  By using
the ideas of Section~\ref{sec:semi-to-uniform}, we are able to pass
the $\chi^2$ test for McL at the 0.05 significance level.  We achieve
$\chi=15.5$ for 20 degrees of freedom using only 15 group operations
per random element.  ($\calR_{15}\inv\overline\calR_{15}$,
$\calR_{15}=g_1^{E_1}\cdots g_{15}^{E_{15}}$, $g_i$ chosen based on
Section~\ref{sec:semi-to-uniform})

% Test(G, partition, preIters, iters, a, b, c, uniformExtra)
% =  Test(G, ConjugacyClasses(G), 15, 1000, 1, 1, 1, 2);

Detailed distributions are provided in the context of the McLaughlin group
in the appendix.

\IGNORE{
{\bf
Base Change paper:  90\% confidence level for covering full group in $x$
iterations.  In the language of statistics, the hypothesis is that the
product replacement algorithm produces a 0.01-uniform random
variable.  We show that for many groups, we can reject that
hypothesis.
}
} % end IGNORE

\IGNORE{

The experimental results are still preliminary.
The groups tested were:
$S_7$ (symmetric group): order 5,040, 2 generators
Fibonacci cube: 181 group operations
product replacement significantly faster

Largest solvable group in $\GL(4,2)$: order 576, 8 generators
62 group operations

$\Syl_2(\GL(4,2))$: order 64, 6 generators -- original gen's random

$\Syl_2(S_8)$:  order 128, 3 generators
148 group operations
product replacement faster

$(Z_2)^7$:  128, 7 generators (generating set already yields uniform
generation)

$S_9$: 362,880, 2 generators, 367 group operations
$\Syl_2(S_4)$: 32,768, 4 generators, 250 group operations

{\bf Construct $g_1,\ldots$, and then look at minimum
and maximum probability density.  Also, report timings
to find $g_1,\ldots$, and timing to take new random group element
based on it.}

} % end IGNORE

\IGNORE{

\section{Erd\H{o}s-R\'enyi generators (constant factor speedup)}

Erd\H os and R\'enyi~\cite{ErdosRenyi65} {\bf did XXX}
Can we get exponential reliability if we use $2\log|G|$ generators?

{\bf
Erd\H os and R\'enyi~\cite[Theorem~1]{ErdosRenyi65} proved that for
\[
   k\ge 2\log|G| +2\log(1/\varepsilon) + \log(1/\delta),
\]
a sequence of $k$ random elements of $G$ will be a sequence of 
$\varepsilon$-uniform Erd\H{o}s - R\'enyi generators with probability
$\ge 1-\delta$.
}

{\bf We can now take $k$ elements from the distribution of $X$.
If one of them is uniform, then the product of the $k$ elements is
uniform by Theorem~\ref{thm:uniform}.  By Chernoff's bound, it's
exponentially likely that we have Erd\H{o}s-R\'enyi generators in
$\calR_j$ as $j$ grows, and if we take a product of $k$ elements for a
random generator, then it is
exponentially likely that it is uniform random as $k$ grows.
Hence, we can use the Erd\H{o}s-R\'enyi theorem directly, with no
assumption about almost random.
}

} % end IGNORE

\IGNORE{

{\bf ALSO TRY:
Assume $U_A$ exists.
\begin{enumerate}
\item If $A\inv A=G$, then $\Pr(A\inv A=g)\ge 1/|G|$.
\item Else $A\inv A\not=G$
\begin{enumerate}
\item $(Y\inv Y)^2$ almost always in $A\inv A$, which should imply
  $A\inv A$ a group (HOW?);  then $r$ a r.s.p. and $A\inv Ar\cap A\inv
  A=\emptyset\Leftarrow A4\cap A=\emptyset$.
\item Can find $g\in (Y\inv Y)^2\setminus A\inv A\Leftarrow |Ag\cup
  A|=2|A|$
\end{enumerate}
\end{enumerate}
}

\subsection{Timing}

{\bf
$O(\log|G|)$ elements in sequence $\calR_1, \calR_2,\ldots$.  Since
$\calR_j=g_1^{E_1} g_2^{E_2} \cdots g_j^{E_j}$ and $E_i\in\{0,1\}$ are
uniform independent random variables, it costs at most $j$ group
operations to find an element from the distribution of $\calR_j$.
Note $j=O(\log|G|)$.  ($j=c\log|G|$ for $c$ depending on probability
of getting the ``right'' $g_i$ at each step.)  The construction of
$g_i$ costs at most $O(\log|G|)$ generators.  (Can put a number on
it.)

So, full timing is $j\times\log(1/\varepsilon)\times$ (cost of
constructing $g_i$) to construct the $g_i$ and
$\log(1/\varepsilon)j\log|G|$ to construct the Erd\H{o}s-R\'enyi
generators, where $\varepsilon$ is the probability of error.
This timing is $O(\log^2|G|\log(1/\varepsilon))$.
}

} % end IGNORE

\section{Product Replacement}
\label{sec:productReplacement}

The Fibonacci cube algorithm can emulate a variation of the product
replacement algorithm, which produces an $\varepsilon$-uniform random
element in $O(\log^2|G|)$ steps.  This should be compared with the
work of Pak~\cite{Pak00} to produce nearly random {\em $k$-sets} (not
elements) in the limiting distribution in $O(\log^9|G|)$ with
$k=O(\log|G|\log\log|G|)$.

To see this, choose $k=O(\log|G|)$ and modify
the product replacement algorithm so that at each step, a randomly
chosen element, $g_i$, of the $k$-set is chosen and {\em all other} elements
of the $k$-set are multiplied by~$g_i$.  Further, after the $i$-th
element has been chosen, it should not be chosen again.  After
$O(\log|G|)$ steps, an element of the $k$-set that has not yet been
chosen will have an $\varepsilon$-uniform distribution.  The proof is
modelled on the proof for the Fibonacci cube algorithm.
Further details will be provided in a different paper.

\section{Permutation Group Membership}
\label{sec:permGroupMembership}

\IGNORE{
{\bf exponentially likely for large base case}
{\bf Do we use any group membership data structure (Knuth, Sims,
Jerrum, deep sift), or can we show something better about deep sift
($D\inv D$) since $D$ is a uniform distribution that we apply
directly?}
} % end IGNORE

Precomputation of a group membership data structure for permutation
groups allows one to compute group orders, find random elements, test
an arbitrary permutation for group membership, etc.  There are at
least four such group membership data structures: Sims's Schreier
vectors (or Schreier trees)~\cite{Sims71}, Knuth's data
structure~\cite{Knuth91}, Jerrum's labelled
branchings~\cite{Jerrum86}, and the deep sift data structure of
Cooperman and Finkelstein~\cite{CoopermanFinkelstein93}.

If $n$ is the permutation degree and $b<n$ is the size of a base, then
$b\le\log|G|\le b\log n$.  Schreier vectors require $O(bn)$ group
operations in the worst case, but $O(\log|G|)$ operations typically,
to produce a random element.  Knuth's data structure and deep sift
require $O(\log|G|)$ operations to produce a random element.  Jerrum's
data structure requires $O(b)$ operations to produce a random element.
While the first three data structures require space proportional to
the time to produce a random element, Jerrum's data structure has the
disadvantage of requiring space for $\theta(n)$ group elements.

Cooperman and
Finkelstein~\cite[Theorem~A]{CoopermanFinkelstein94} had previously
demonstrated a random base change algorithm requiring $O(\log|G|)$
random group elements as input.  The base change algorithm produces a
group membership data structure, thus solving permutation group
membership.  The original paper assumed that the random group elements
came from a Schreier vector, but the result of this paper provides an
alternative source of such random elements.  Combining this paper with
the random base change algorithm and any of the four group membership
data structures yields a Monte Carlo group membership algorithm operating in
$O(\log^2|G|)$ group operations.

\begin{sloppypar}
Prior to this, the fastest general algorithm was the deep sift
algorithm of Cooperman and Finkelstein, requiring $O(n^2\log^4n)$
group operations, and the fastest small base algorithm Babai,
Cooperman, Finkelstein and Seress~\cite{BabaiCoopermanEtal91b},
required $O(\log^3|G|+b^2\log^2|G|\log(b+\log n)|b^2(\log
b)(\log^3|G|)(\log n)/n)$ group operations.  In both cases, the
$O(\log^2|G|)$ group operations using the new Fibonacci cube algorithm
represents a significant improvement.
\end{sloppypar}

Any of these Monte Carlo algorithms can be upgraded to Las Vegas by
applying a strong generating test afterwards.  The danger with Monte
Carlo algorithms is that they may not produce enough group elements to
form a full strong generating set.  Cooperman and
Finkelstein~\cite{CoopermanFinkelstein91} demonstrate a $O(\log^3|G|)$
algorithm for testing if a set of group elements forms a strong
generating set.  (In fact, the algorithm is $O(n^4)$ for a permutation
group acting on $n$~points.)

% end IGNORE

\section{Conclusion}

The Fibonacci cube algorithm has been demonstrated to produce a
$\gamma$-uniform random variable in $O(\log^2|G|)$ group operations.
From that distribution $\varepsilon$-uniform elements with
$O((\log1/\varepsilon)\log|G|)$ group operations can be computed.  The
algorithm is asymptotically faster than previous theoretical
algorithms and also empirically faster than the product replacement
heuristic for many groups.  The faster random generation algorithm
also yields a faster permutation group membership algorithm.

The coefficient of complexity of the Fibonacci cube algorithm analyzed
in this paper is still unacceptably high.  This may not be an issue
for computations that have an independent check for correctness, such
as Las Vegas algorithms, since the experimental results are
competitive.  An expanded version of this paper will refine the
analysis to produce a smaller coefficient of complexity.

The large coefficient arises due to the constant 0.997 arising from
Section~\ref{sec:escape}.  In Theorem~\ref{thm:escape}, we are too
greedy in demanding that $A'A'=A'={A'}\inv$ (and therefore
$A'g\setminus A'=\emptyset$).  If we prove only that $A'A'$ does not
differ greatly from $A'$, we can still prove that a random subproduct
has reasonable probability of allowing us to escape the set~$A'$.

\section{Acknowledgement}
The author thanks Igor Pak for reading an early version of this
manuscript.

\IGNORE{
{\bf The coefficient of complexity is {\bf XXX}.  The large
  coefficient of complexity arises from the small constant 0.006 in
  Remark~\ref{rem:constants}.  This occurs because we are too greedy
  in demanding that $A'A'=A'={A'}\inv$ in Theorem~\ref{thm:escape}
  (and therefore $A'g\setminus A'=\emptyset$).  If we prove only that
  $A'A'$ does not differ greatly from $A'$, then for $g$ a random
  subproduct on the generators, $|A'g\setminus A'|/|A'|$ will be
  sufficiently large to be used in the algorithms of the rest of the
  paper.  This will be addressed in a future analysis.
}
} % end IGNORE

\bibliographystyle{alpha}
\bibliography{rand-grp-elts}

\section*{Appendix: Computational Experiment}

This appendix is a quick note on a computation suggested by Persi
Diaconis.  It is about a quick computational experiment.
It is not intended to be a polished document.

I test McLaughlin's group (McL).  I compare the
true distribution of elements according to the conjugacy classes,
with the distribution according to conjugacy classes produced
by the random generator of the paper.

I take the constants $a=b=c=1$ in the Fibonacci Cube algorithm.
I use only the Fibonacci Cube
algorithm, which in principle produces only a semi-uniform
random variable.
That is, in principle, this distribution will satisfy only
\[
 \forall g \in G,\quad
  \Pr(X = g) > \alpha/|G|
\]

The paper has an additional step for producing nearly uniform
random variables.  I will test the full algorithm at a later date.
I suspect the full algorithm will represent an improvement.
But for now even the semi-uniform random variables seem
to be close enough to uniform.

The code was written using GAP~4.2.
The test here is for
McL (McLaughlin group), of order 898,128,000, with 2~generators,
        based on a permutation representation on 275~points.
The representation is
        provided by Walter Kim, U. Chicago, Feb., 2000.

For McLaughlin's group, there are 24 conjugacy classes.
For each conjugacy class, $C_i$, I compute an integer,
$\lceil C_i/10^6 \rceil$.  This is for convenience, since
GAP doesn't handle floating point.  Since $\sum \lceil
C_i/10^6\rceil=886$, I test
the random generator by generating exactly 886 elements, and test
their distribution into conjugacy classes.

In each case, the first row is the distribution of elements produced by
the random generator (the number of elements in each of the 24 conjugacy
classes).  The second row corresponds to the true distribution,
normalized to the form $\lfloor C_i/10^6 \rfloor$.
The notation 30~terms means that $\calR_{30}$ was computed in the
notation of the paper.  The $O(\log^2|G|)$ precomputation refers to the
computation of $g_1,\ldots,g_{30}$ for $\calR_{30}=g_1^{E_1}\cdots
g_{30}^{E_{30}}$.  The 886 random elements are then each drawn
from $\calR_{30}\inv\overline\calR_{30}$.  This is the $O(\log|G|)$
computation.  On average, the $O(\log|G|)$ computation of a random
element from $\calR_{30}\inv\overline\calR_{30}$ costs 30~group
operations per random element (29~multiplications and one inverse).

Note that for less than 20 terms,
there are too many pseudo-random elements in the first,
third, tenth and eleventh conjugacy classes.  This experimental
observation reflects the theoretical model, which states only that
$\calR_i\inv\overline\calR_i$ is semi-uniform.  A future experiment
will also test the theory of Section~\ref{sec:semi-to-uniform} for
converting semi-uniform to $\varepsilon$-uniform.  This should be more
efficient in producing $\varepsilon$-uniform random elements.

{\small
\begin{verbatim}
Experiment 1:
  30 terms, 174 group operations for O(log^2|G|) precomputation.  30 ops/rand elt.
[ 0, 38, 0, 31, 20, 32, 22, 0, 43, 0, 2, 12, 82, 75, 55, 68, 62, 72, 78, 102, 2, 31, 25, 34 ] 
[ 0, 35, 1, 29, 29, 29, 29, 0, 29, 0, 2, 9, 74, 64, 64, 64, 64, 81, 81, 112, 0, 33, 33, 24 ]

Experiment 2:
  30 terms, 144 group operations for O(log^2|G|) precomputation.  30 ops/rand elt.
[ 0, 31, 0, 33, 34, 27, 19, 0, 24, 0, 1, 13, 76, 79, 72, 64, 59, 94, 75, 105, 3, 29, 28, 20 ]
[ 0, 35, 1, 29, 29, 29, 29, 0, 29, 0, 2, 9, 74, 64, 64, 64, 64, 81, 81, 112, 0, 33, 33, 24 ]

Experiment 3:
  20 terms, 74 group operations for O(log^2|G|) precomputation.  20 ops/rand elt.
[ 0, 46, 3, 25, 23, 27, 35, 0, 19, 0, 2, 7, 66, 69, 49, 50, 71, 94, 72, 137, 1, 41, 29, 20 ]
[ 0, 35, 1, 29, 29, 29, 29, 0, 29, 0, 2, 9, 74, 64, 64, 64, 64, 81, 81, 112, 0, 33, 33, 24 ]

Experiment 4:
  20 terms, 88 group operations for O(log^2|G|) precomputation.  20 ops/rand elt.
[ 0, 32, 1, 32, 31, 31, 18, 0, 22, 1, 3, 11, 78, 66, 49, 67, 76, 81, 74, 120, 0, 39, 23, 31 ]
[ 0, 35, 1, 29, 29, 29, 29, 0, 29, 0, 2, 9, 74, 64, 64, 64, 64, 81, 81, 112, 0, 33, 33, 24 ]

Experiment 5:
  15 terms, 48 group operations for O(log^2|G|) precomputation.  15 ops/rand elt.
[ 1, 44, 23, 38, 32, 27, 33, 0, 25, 9, 8, 11, 72, 50, 59, 50, 56, 79, 72, 104, 4, 31, 33, 25 ]
[ 0, 35, 1, 29, 29, 29, 29, 0, 29, 0, 2, 9, 74, 64, 64, 64, 64, 81, 81, 112, 0, 33, 33, 24 ]

Experiment 6:
  15 terms, 44 group operations for O(log^2|G|) precomputation.  15 ops/rand elt.
[ 8, 41, 54, 44, 36, 12, 14, 0, 10, 17, 11, 6, 94, 34, 46, 50, 69, 82, 71, 54, 1, 38, 52, 42 ]
[ 0, 35, 1, 29, 29, 29, 29, 0, 29, 0, 2, 9, 74, 64, 64, 64, 64, 81, 81, 112, 0, 33, 33, 24 ]

Experiment 7:
  10 terms, 21 group operations for O(log^2|G|) precomputation.  10 ops/rand elt.
[ 0, 47, 54, 55, 53, 11, 14, 1, 6, 30, 8, 11, 96, 47, 66, 44, 42, 80, 74, 39, 4, 45, 37, 22 ]
[ 0, 35, 1, 29, 29, 29, 29, 0, 29, 0, 2, 9, 74, 64, 64, 64, 64, 81, 81, 112, 0, 33, 33, 24 ]

Experiment 8:
  10 terms, 19 group operations for O(log^2|G|) precomputation.  10 ops/rand elt.
[ 17, 49, 82, 27, 28, 10, 12, 0, 15, 30, 7, 7, 101, 36, 26, 40, 54, 110, 107, 12, 0, 46, 51, 19 ]
[ 0, 35, 1, 29, 29, 29, 29, 0, 29, 0, 2, 9, 74, 64, 64, 64, 64, 81, 81, 112, 0, 33, 33, 24 ]

\end{verbatim}
}
\end{document}

%% file: Am.pstex_t
\begin{picture}(0,0)%
\includegraphics{Am.pstex}%
\end{picture}%
\setlength{\unitlength}{3158sp}%
\begingroup\makeatletter\ifx\SetFigFont\undefined%
\gdef\SetFigFont#1#2#3#4#5{%
  \reset@font\fontsize{#1}{#2pt}%
  \fontfamily{#3}\fontseries{#4}\fontshape{#5}%
  \selectfont}%
\fi\endgroup%
\begin{picture}(5283,2024)(1351,-3548)
\end{picture}

%% file: Am-case1a.pstex_t
\begin{picture}(0,0)%
\special{psfile=Am-case1a.pstex}%
\end{picture}%
\setlength{\unitlength}{2763sp}%
\begingroup\makeatletter\ifx\SetFigFont\undefined
% extract first six characters in \fmtname
\def\x#1#2#3#4#5#6#7\relax{\def\x{#1#2#3#4#5#6}}%
\expandafter\x\fmtname xxxxxx\relax \def\y{splain}%
\ifx\x\y   % LaTeX or SliTeX?
\gdef\SetFigFont#1#2#3{%
  \ifnum #1<17\tiny\else \ifnum #1<20\small\else
  \ifnum #1<24\normalsize\else \ifnum #1<29\large\else
  \ifnum #1<34\Large\else \ifnum #1<41\LARGE\else
     \huge\fi\fi\fi\fi\fi\fi
  \csname #3\endcsname}%
\else
\gdef\SetFigFont#1#2#3{\begingroup
  \count@#1\relax \ifnum 25<\count@\count@25\fi
  \def\x{\endgroup\@setsize\SetFigFont{#2pt}}%
  \expandafter\x
    \csname \romannumeral\the\count@ pt\expandafter\endcsname
    \csname @\romannumeral\the\count@ pt\endcsname
  \csname #3\endcsname}%
\fi
\fi\endgroup
\begin{picture}(5111,1812)(1351,-3486)
\end{picture}

%% file: Am-case1b.pstex_t
\begin{picture}(0,0)%
\includegraphics{Am-case1b.pstex}%
\end{picture}%
\setlength{\unitlength}{2763sp}%
\begingroup\makeatletter\ifx\SetFigFont\undefined%
\gdef\SetFigFont#1#2#3#4#5{%
  \reset@font\fontsize{#1}{#2pt}%
  \fontfamily{#3}\fontseries{#4}\fontshape{#5}%
  \selectfont}%
\fi\endgroup%
\begin{picture}(5283,1812)(1351,-3486)
\end{picture}

%% file: Am-case2a.pstex_t
\begin{picture}(0,0)%
\includegraphics{Am-case2a.pstex}%
\end{picture}%
\setlength{\unitlength}{2763sp}%
\begingroup\makeatletter\ifx\SetFigFont\undefined%
\gdef\SetFigFont#1#2#3#4#5{%
  \reset@font\fontsize{#1}{#2pt}%
  \fontfamily{#3}\fontseries{#4}\fontshape{#5}%
  \selectfont}%
\fi\endgroup%
\begin{picture}(5283,1812)(1351,-3486)
\end{picture}

%% file: Am-case2b.pstex_t
\begin{picture}(0,0)%
\includegraphics{Am-case2b.pstex}%
\end{picture}%
\setlength{\unitlength}{2763sp}%
\begingroup\makeatletter\ifx\SetFigFont\undefined%
\gdef\SetFigFont#1#2#3#4#5{%
  \reset@font\fontsize{#1}{#2pt}%
  \fontfamily{#3}\fontseries{#4}\fontshape{#5}%
  \selectfont}%
\fi\endgroup%
\begin{picture}(5283,1829)(1351,-3503)
\end{picture}

%% file: Am-case3.pstex_t
\begin{picture}(0,0)%
\includegraphics{Am-case3.pstex}%
\end{picture}%
\setlength{\unitlength}{3158sp}%
\begingroup\makeatletter\ifx\SetFigFont\undefined%
\gdef\SetFigFont#1#2#3#4#5{%
  \reset@font\fontsize{#1}{#2pt}%
  \fontfamily{#3}\fontseries{#4}\fontshape{#5}%
  \selectfont}%
\fi\endgroup%
\begin{picture}(5283,1812)(1351,-3486)
\end{picture}

%% file: rand-grp-elts.bbl
\newcommand{\etalchar}[1]{$^{#1}$}
\begin{thebibliography}{CLGM{\etalchar{+}}95}

\bibitem[Bab91]{Babai91}
L.~Babai.
\newblock Local expansion of vertex-transitive graphs and random generation in
  finite groups.
\newblock In {\em 23$^{\mbox{rd}}$ ACM Symposium on Theory of Computing
  (STOC)}, pages 164--174, New York, 1991. Association for Computing Machinery.

\bibitem[Bab95]{Babai95}
L.~Babai.
\newblock Automorphism groups, isomorphism, reconstruction.
\newblock In {\em Handbook of Combinatorics}, volume~2, pages 1447--1540.
  Elsevier, Amsterdam, 1995.

\bibitem[BCF{\etalchar{+}}91]{BabaiCoopermanetal91}
L.~Babai, G.~Cooperman, L.~Finkelstein, E.~M.\ Luks, and \'Akos Seress.
\newblock Fast {M}onte {C}arlo algorithms for permutation groups.
\newblock In {\em Proc.\ $23$rd ACM STOC}, pages 90--100, 1991.

\bibitem[BCFS91]{BabaiCoopermanEtal91b}
L.~Babai, G.~Cooperman, L.~Finkelstein, and \'Akos Seress.
\newblock Nearly linear time algorithms for permutation groups with a small
  base.
\newblock In {\em Proc. of International Symposium on Symbolic and Algebraic
  Computation {\rm ISSAC '91}}, pages 200--209. (Bonn), ACM Press, 1991.

\bibitem[BP02]{BabaiPak02}
L.~Babai and I.~Pak.
\newblock Strong bias of group generators: an obstacle to the ``product
  replacement algorithm''.
\newblock {\em J.~Algorithms}, 2002.
\newblock to appear.

\bibitem[BS84]{BabaiSzemeredi84}
L.~Babai and E.~Szemer\'edi.
\newblock On the complexity of matrix group problems, {I}.
\newblock In {\em Proc. $25$th IEEE Symposium on Foundations of Computer
  Science}, pages 229--240, 1984.

\bibitem[CF91]{CoopermanFinkelstein91}
G.~Cooperman and L.~Finkelstein.
\newblock A strong generating test and short presentations for permutation
  groups.
\newblock {\em J.~of Symbolic Computation}, 12:475--497, 1991.

\bibitem[CF93]{CoopermanFinkelstein93}
G.~Cooperman and L.~Finkelstein.
\newblock Combinatorial tools for computational group theory.
\newblock In {\em Groups and Computation}, volume~11 of {\em Amer.\ Math.\
  Soc.\ DIMACS Series}, pages 53--86. (DIMACS, 1991), 1993.

\bibitem[CF94]{CoopermanFinkelstein94}
G.~Cooperman and L.~Finkelstein.
\newblock A random base change algorithm for permutation groups.
\newblock {\em J.~of Symbolic Computation}, 17:513--528, 1994.

\bibitem[Che52]{Chernoff52}
H.~Chernoff.
\newblock A measure of asymptotic efficiency for tests of a hypothesis based on
  the sum of observations.
\newblock {\em Annals of Math. Statistics}, 23:493--507, 1952.

\bibitem[CLGM{\etalchar{+}}95]{CellerLeedhamGreenEtAl95}
F.~Celler, C.R. Leedham-Green, S.~Murray, A.~Niemeyer, and E.A. O'Brien.
\newblock Generating random elements of a finite group.
\newblock {\em Communications of Algebra}, 23:4931--4948, 1995.

\bibitem[DSC93]{DiaconisSaloffCoste93}
P.~Diaconis and L.~Saloff-Coste.
\newblock Comparison techniques for random walk on finite groups.
\newblock {\em Ann. Prob.}, 21:2131--2156, 1993.

\bibitem[DSC98]{DiaconisSaloffCoste98}
P.~Diaconis and L.~Saloff-Coste.
\newblock Walks on generating sets of groups.
\newblock {\em Inventiones Mathematicae}, 134:251--299, 1998.

\bibitem[Jer86]{Jerrum86}
M.~Jerrum.
\newblock A compact representation for permutation groups.
\newblock {\em J.~Algorithms}, 7:60--78, 1986.

\bibitem[Knu91]{Knuth91}
D.E. Knuth.
\newblock Notes on efficient representation of perm groups.
\newblock {\em Combinatorica}, 11:57--68, 1991.
\newblock preliminary version since 1981.

\bibitem[LG01]{Leedham-Green01}
C.R. Leedham-Green.
\newblock The computational matrix group project.
\newblock In {\em Groups and Computation~III}, pages 123--138. DeGruyter
  Publishers, 2001.

\bibitem[Pak00]{Pak00}
I.~Pak.
\newblock The product replacement algorithm is polynomial.
\newblock In {\em Proc.~41$^{\mbox{st}}$ {IEEE} Symposium on Foundations of
  Computer Science (FOCS)}, pages 476--485. IEEE Press, 2000.

\bibitem[Sim71]{Sims71}
C.C. Sims.
\newblock Computation with permutation groups.
\newblock In {\em Proc. Second Symposium on Symbolic and Algebraic
  Manipulation}, pages 23--28. ACM Press, 1971.

\end{thebibliography}
